\documentclass[MSNbibl,number,citesort,seceqn,dvips]{arxbj}
\usepackage{upgreek}
\usepackage{graphicx}

\aid{0}
\volume{20}
\issue{2}
\pubyear{2014}
\firstpage{803}
\lastpage{845}
\doi{10.3150/13-BEJ507} 

\makeatletter

\newcommand{\rrvert}{\vert}
\newcommand{\llvert}{\vert}
\newcommand{\eqref}[1]{(\ref{#1})}
\newcommand{\ex}{\mathrm{e}}
\newtheorem{lemma}{Lemma}[section]
\newtheorem{theorem}[lemma]{Theorem}
\newcommand{\bbc}{\mathbb{C}}
\newtheorem{proposition}[lemma]{Proposition}
\newtheorem{corollary}[lemma]{Corollary}
\newremark{remark}[lemma]{Remark}

\newcommand{\var}{\operatorname{var}}
\newcommand{\cov}{\operatorname{cov}}
\newcommand{\corr}{\operatorname{corr}}

\newcommand{\stp}{\stackrel{P}{\rightarrow}}
\newcommand{\std}{\stackrel{d}{\rightarrow}}
\newcommand{\stv}{\stackrel{v}{\rightarrow}}
\newcommand{\eqd}{\stackrel{d}{=}}

\newcommand{\w}{\omega}
\newcommand{\nto}{n\to\infty}
\newcommand{\kto}{k\to\infty}
\newcommand{\xto}{x\to\infty}
\newcommand{\ov}{\overline}
\newcommand{\wt}{\widetilde}
\newcommand{\wh}{\widehat}
\newcommand{\vep}{\varepsilon}
\newcommand{\la}{\lambda}

\newcommand{\bbr}{\mathbb{R}}
\newcommand{\bbz}{\mathbb{Z}}

\newcommand{\bfZ}{\mathbf{Z}}
\newcommand{\bfc}{\mathbf{c}}
\def\mid{\vert}
\makeatother

\begin{document}
\begin{frontmatter}

\title{A Fourier analysis of extreme events}

\runtitle{A Fourier analysis of extreme events}

\begin{aug}
\author{\fnms{Thomas} \snm{Mikosch}\thanksref{e1}\ead[label=e1,mark]{mikosch@math.ku.dk}}%
\and
\author{\fnms{Yuwei} \snm{Zhao}\corref{}\thanksref{e2}\ead[label=e2,mark]{zhao@math.ku.dk}}
\runauthor{T. Mikosch and Y. Zhao} 
\address{University of Copenhagen, Department of Mathematics,
Universitetsparken 5,
DK-2100 Copenhagen, Denmark. \printead{e1,e2}}
\end{aug}

\received{\smonth{4} \syear{2012}}
\revised{\smonth{11} \syear{2012}}

%
\begin{abstract}
The \emph{extremogram} is an asymptotic correlogram for extreme
events constructed
from a regularly varying stationary sequence.
In this paper, we define a frequency domain analog of the correlogram: a
periodogram generated from a suitable sequence of indicator functions
of rare events.
We derive basic properties of the periodogram such as the asymptotic independence at the Fourier frequencies and use this property to
show that weighted versions of the periodogram are consistent
estimators of
a spectral density derived from the extremogram.
\end{abstract}

%

\end{frontmatter}

\section{Introduction}\label{secintrod}
In this paper, we study an analog of the periodogram for extremal events.
In classical time series analysis, the periodogram is a method of
moments estimator of the
spectral density of a second order stationary time series $(X_t)$;
see, for example, the standard monographs
Brillinger~\cite{brillinger1981},
Brockwell and Davis \cite{brockwelldavis1991},
Grenander and Rosenblatt \cite{grenanderrosenblatt1984},
Hannan \cite{hannan1960},
Priestley \cite{priestley1981}. The notions of spectral
density and periodogram are the respective frequency domain analogs of
the autocorrelation function
and the
sample autocorrelation function in the time domain. In the context of
extremal events,
these notions are not meaningful since second order characteristics
are not suited for describing the occurrence of rare events.

However, Davis and Mikosch \cite{davismikosch2009} introduced
a time domain analog of the autocorrelation function,
the \emph{extremogram} for rare events. For an $\mathbb{R}^{d}$-valued
strictly stationary
time series $(X_t)$ and a Borel set $A$ bounded away from zero,
\emph{the extremogram at lag $h\ge0$} is given as the limit
%
\begin{equation}
\label{eq1} \rho_{A}(h)=\lim_{\xto} P
\bigl(x^{-1} X_h\in A\mid x^{-1} X_0
\in A\bigr).
\end{equation}
This definition requires that the support of $X$ (here and in what
follows, $X$ denotes a generic element of any stationary sequence $(X_t)$)
is unbounded and, more importantly, that the limit on the right-hand
side exists. In general, these limits do not exist. A~sufficient condition
for their existence is \emph{regular variation} of all pairs $(X_0,X_h)$
or, more
generally,
\emph{regular variation of the
finite-dimensional distributions} of the process $(X_t)$. A~precise
definition of regular variation will
be given in Section~\ref{subsecregvar}. Since $A$ is assumed to be
bounded away from zero, the probabilities $P(x^{-1} X\in A)$ converge to
zero as $\xto$. Then the following calculation is straightforward for~$A$:
\begin{eqnarray*}
\lim_{\xto} \corr(I_{\{x^{-1} X_0\in A\}},I_{\{x^{-1} X_h\in A\}
})&=&\lim
_{\xto}\frac{
P(x^{-1} X_0\in A ,x^{-1} X_h\in A)-[P(x^{-1} X\in A)]^2}{P(x^{-1}
X\in A)(1- P(x^{-1}
X\in A))}
\\
&=&\lim_{\xto} P\bigl(x^{-1} X_h\in A
\mid x^{-1} X_0\in A\bigr)=\rho_{A}(h).
\end{eqnarray*}
For fixed $x$, $(I_{\{x^{-1} X_t\in A\}})_{t\in\bbz}$ constitutes a strictly
stationary sequence. The limit sequence $(\rho_{A}(h))$ inherits
the property of correlation function from $(\corr(I_{\{x^{-1} X_0\in
A\}},I_{\{x^{-1} X_h\in A\}}))$.
Therefore, in an asymptotic sense, one can use the notions of classical time series analysis (such
as the autocorrelation function)
for the sequences of indicator functions $(I_{\{x^{-1} X_t\in A\}
})_{t\in
\bbz}$.
Of course, there are several
crucial differences to classical time series analysis.
\begin{itemize}
\item
The notion of autocorrelation function is only
defined in an asymptotic sense.
\item
The strictly stationary sequence of indicator functions $(I_{\{x^{-1}
X_t\in
A\}})_{t\in\bbz}$
depends on the threshold $x$, that is, we are dealing with an array of
strictly stationary processes.
\item
By definition, the values $\rho_{A}(h)$ cannot be negative.
\end{itemize}
Davis and Mikosch \cite{davismikosch2009,davismikosch2011}
introduced the extremogram and calculated the extremogram for various
standard regularly varying time series models
such as the GARCH model, stochastic volatility and linear processes
with regularly varying noise, and infinite variance stable processes; see also
Section~\ref{subsecexam}. They studied the basic asymptotic properties of
the extremogram
(consistency, asymptotic normality) and
also introduced a frequency domain
analog of the correlation function~$\rho_{A}$ given
as the Fourier series
%
\begin{equation}
\label{eq5} f_A(\la)= \sum_{h\in\bbz}
\rho_{A}(h) \ex^{-\mathrm{i} h \la} ,\qquad \la \in[0,\uppi].
\end{equation}
A natural estimator of $
f_A(\la)$ is
found by replacing the correlations $\rho_{A}(h)$ by sample analogs.
The convergence in the mean square sense of such an analog of the
classical periodogram estimator towards the spectral density
$f_A(\la)$ at a fixed frequency $\la$ was shown in \cite{davismikosch2009}.
However, the periodogram of $(I_{\{x^{-1} X_t\in A\}})_{t\in\bbz}$
used in
\cite{davismikosch2009} had to be truncated to achieve
consistency; the truncation level depended on some mixing rate which
is unknown for real-life data.
In this paper, we overcome this inconvenience. In addition, we study
the periodogram ordinates of the indicator functions at finitely many frequencies.
We show
that the limiting vector
of the periodogram ordinates at distinct fixed or Fourier frequencies
converges in distribution to
a vector of independent exponential random variables.
This property parallels the asymptotic theory for the periodogram of a second order stationary sequence; see,
for example, Brockwell and Davis \cite{brockwelldavis1991}, Chapter~10.

In
classical time series analysis, the asymptotic independence of the
periodogram at distinct frequenc
ies is the theoretical
basis for consistent estimation of the spectral density via weighted
averages or kernel based methods. We show that weighted average
estimators of the periodogram evaluated at Fourier frequencies in the
neighborhood of a fixed non-zero frequency are consistent estimators of
the limiting spectral density.

The paper is organized as follows. In Section~\ref{secprelim}, we
introduce basic notions and conditions used throughout this paper.
In Section~\ref{subsecregvar}, we define regular variation of a strictly
stationary sequence. In Section~\ref{subsecmix}, we consider those mixing
conditions which are relevant for the results of this paper.
The periodogram of extreme events is introduced in
Section~\ref{subsecper}.
In Section~\ref{subsecexam}, we discuss some regularly varying strictly
stationary sequences. Among them are linear, stochastic volatility and
max-moving average
processes with regularly varying noise. We give
expressions for the extremogram and, if possible, for the
corresponding spectral density. In Section~\ref{secprop}, we give the
main results of this paper. We start in Section~\ref{subsecperiod}
by showing that the periodogram ordinates of extreme events are asymptotically
uncorrelated
at distinct fixed or Fourier frequencies in the interval
$(0,\uppi)$. Next, in Section~\ref{subsecclt} we show that the
periodogram ordinates at distinct fixed or Fourier frequencies converge to
independent exponential random variables. This property is exploited in
Section~\ref{subsecsmooth} to show that weighted averages of
periodogram ordinates evaluated at Fourier frequencies in a small neighborhood
of a fixed frequency yield consistent estimates of the underlying
spectral density at the given frequency.
In Section~\ref{sec6}, we give a short discussion of work related to
the extremogram or the spectral analysis of sequences of indicator
functions.
The proofs depend on various
calculations involving formulas for sums of trigonometric functions.
Some of these formulas and related calculations are given in the \hyperref[secappendix]{Appendix}.

\section{Preliminaries}\label{secprelim}
\subsection{Regular variation}\label{subsecregvar}
It was mentioned in Section~\ref{secintrod} that one needs
conditions to ensure that the limits $\rho_{A}(h)$
in \eqref{eq1} exist. A sufficient condition for this to
hold is \emph{regular variation} of the strictly stationary sequence $(X_t)$.
Regular variation is a convenient tool for modeling multivariate
heavy-tail phenomena and serial extremal dependence in a time series;
see Resnick's monographs \cite{resnick1987,resnick2007},
Resnick~\cite{resnick1986},
Basrak and Segers
\cite{basraksegers2009,basrakkrizmanicsegers2011}, Davis and
Hsing \cite{davishsing1995}, Embrechts \textit{et al.} \cite
{embrechtskluppelbergmikosch1997},
Jakubowski \cite{jakubowski1993,jakubowski1997},
Bartkiewicz \textit{et
al.}~\cite{bartkiewiczjakubowskimikoschwintenberger2011}, and the
references therein. Regular variation is particularly useful for
modeling extremes in financial time series; see Basrak \textit{et
al.} \cite{basrakdavismikosch2002},
Mikosch and St\u aric\u a~\cite{mikoschstarica2000},
Davis and Mikosch \cite
{davismikosch2001,davismikosch2009a,davismikosch2009b};
cf. Andersen \textit{et al.} \cite{andersendaviskreissmikosch2009} and the
references therein. See also the examples in Section~\ref{subsecexam}.

A random vector $X$ with values in $\bbr^d$ for some $d\ge1$ is
\emph{regularly varying} if there exists a non-null Radon measure $\mu
$ on the
Borel $\sigma$-field of $ \ov\bbr^d_0=\ov\bbr^d\setminus\{\mathbf{0}\}$,
where $\ov\bbr=\bbr\cup\{\infty,-\infty\}$, such that
%
\begin{equation}
\label{eq2} \frac{P(x^{-1}X\in\cdot)}{P(|X|>x)}\stv\mu(\cdot) ,\qquad \xto.
\end{equation}
Here $\stv$ denotes \emph{vague convergence} on the Borel $\sigma
$-field of $
\ov\bbr^d_0$; for definitions see Kallenberg~\cite{kallenberg1983},
Resnick \cite{resnick1987,resnick1986}.
In this context, bounded
sets are those which are bounded away from zero and the Radon measure $\mu$
charges finite mass to these sets. Then, necessarily, there exists an
$\alpha\ge0$ such that $\mu(t A)= t^{-\alpha}\mu(A)$, $t>0$, for
all $A$
in the Borel $\sigma$-field of $\ov\bbr^d_0$. We refer to \emph
{regular variation of $X$ with limiting measure $\mu$ and index $\alpha$.}
A multivariate $t$-distributed random vector is regularly varying and
the index $\alpha$ is the degree of freedom. Other well known multivariate
regularly varying distributions are the multivariate $F$- and Fr\'echet distributions; see Resnick
\cite{resnick1987}, Chapter 5, in particular Section 5.4.2.

We will often use
an equivalent sequential version of \eqref{eq2}:
there exists $(a_n)$ such that $a_n\to\infty$ as $\nto$ and
%
\begin{equation}
\label{eq3} n P\bigl(a_n^{-1} X\in\cdot\bigr)\stv\mu(\cdot)
, \qquad\nto.
\end{equation}
A possible choice of $(a_n)$ is given by the $(1-1/n)$-quantile of
$|X|$.

Now, a strictly stationary $d$-dimensional sequence $(X_t)$ is \emph
{regularly varying} if the lagged vectors $Y_h= \mathrm{vec}(X_0,
\ldots,
X_{h})$, $h\ge0$, are
regularly varying with index $\alpha$. Of course, the limiting
non-null Radon
measures $\mu_h$ in
\eqref{eq2} now depend on the lag $h$ and the normalization in
\eqref{eq3} would also change with $h$. In the context of
this paper it is convenient to choose the
normalizations of the rare event probabilities
independently of $h$. In particular, we will use the following relations
for $h\ge0$,
\begin{eqnarray*}
\frac{P(x^{-1}Y_h\in\cdot)}{P(|X_0|>x)}&\stv& \mu_h(\cdot) ,\qquad \xto,
\\
n P\bigl(a_n^{-1} Y_h\in\cdot\bigr)&\stv&
\mu_h(\cdot) ,\qquad \nto,
\end{eqnarray*}
where $(a_n)$ satisfies $n P(|X_0|>a_n)\to1$, as $\nto$. These
relations are equivalent to the definitions \eqref{eq2} and
\eqref{eq3} of regular variation of $Y_h$.

Now we are in the position to verify that the limits $\rho_{A}(h)$ in
\eqref{eq1} exist for any Borel set $A\subset\ov\bbr^d_0$ bounded
away from zero. Write $\wt A= A \times\ov
\bbr_0^{dh}$ and $\wt B=A\times\ov\bbr_0^{d(h-1)} \times
A$. These sets are bounded away from zero in $\ov\bbr^{d(h+1)}_0$. If
these sets are continuity sets with respect to $\mu_{h}$ we obtain
from the
sequential definition of regular variation of $Y_h$ for $h\ge0$,
\begin{eqnarray*}
\rho_{A}(h) &=& \lim_{n\to\infty} P\bigl(a_n^{-1}
X_h \in A \mid a_n^{-1} X_0 \in A
\bigr)
\\
&=& \lim_{n\to\infty} \frac{n P(a_n^{-1} Y_h \in\wt B)}{n P(a_n^{-1}
Y_h \in\wt A)}= \frac{\mu_h (\wt B)}{\mu_h (\wt A)}.
\end{eqnarray*}
%
\subsection{The mixing and dependence conditions \textup{(M)}, \textup{(M1)}
and \textup{(M2)}}\label{subsecmix}
The results in Davis and Mikosch
\cite{davismikosch2009,davismikosch2011}
were proved under the following mixing/dependence condition on the
sequence $(X_t)$.

\begin{longlist}[(M)]
\item[(M)] The sequence $(X_t)$ is strongly mixing with rate
function $(\xi_t)$. There exist $m=m_n\to\infty$ and
$r_n\to\infty$
such that $m_n/n\to0$ and $r_n/m_n\to0$ and
%
\begin{equation}
\label{eqll} \lim_{\nto} m_n\sum
_{h=r_n}^\infty\xi_h=0 ,
\end{equation}
and for all $\epsilon>0$,
%
\begin{eqnarray}
\label{eq33} \lim_{k\to\infty}\limsup_{\nto}
m_n\sum_{h=k}^{r_n}
P\bigl(|X_h|>\epsilon a_m ,|X_0|>\epsilon
a_m\bigr)=0.
\end{eqnarray}
\end{longlist}
Condition \eqref{eq33} is similar in spirit to condition (2.8)
used in Davis and Hsing \cite{davishsing1995} for establishing
convergence of a sequence of point processes to a limiting cluster
point process. It is much
weaker than the anti-clustering condition $D'(\epsilon a_n)$ of
Leadbetter
which is well known in the extreme value literature; see Leadbetter \textit{et
al.} \cite{leadbetterlindgrenrootzen1983} or Embrechts \textit{et al.}
\cite{embrechtskluppelbergmikosch1997}. Since we choose
$(a_n)$ such that $n P(|X|>a_n)\to1$ as $\nto$, \eqref{eq33} is
equivalent to
\[
\lim_{k\to\infty}\limsup_{\nto} \sum
_{h=k}^\infty P\bigl(|X_h|>\epsilon
a_m\mid|X_0|>\epsilon a_m\bigr)=0 ,\qquad
\epsilon>0.
\]
In addition, we also need the following technical condition, using the same
notation as in (M).
\begin{longlist}[(M1)]
\item[(M1)]
The sequences $(m_n)$, $(r_n)$, $k_n=[n/m_n]$ from (M) also satisfy
the growth conditions $k_n \xi_{r_n} \to0$, and $m_n=\mathrm{o}(n^{1/3})$.
\end{longlist}
%
\begin{remark}\label{remrm1}
Some of the examples in Section~\ref{subsecexam} are strongly mixing
with geometric rate, that is, there exists $a\in(0,1)$ such that $\xi
_h\le
a^h$ for sufficiently large $h$.
Then \eqref{eqll} is satisfied if $m_n a^{r_n}=\mathrm{o}(1)$.
If $m_n=n^\gamma$ for
some $\gamma\in(0,1)$ then \eqref{eqll} is
satisfied for $r_n=c\log n$ if $c$ is chosen sufficiently large and
(M1) trivially holds as well. If $\xi_h\le h^{-s}$ for some $s>1$ and
sufficiently large $h$ then \eqref{eqll} is satisfied if
$m_n r_n^{-s+1}=\mathrm{o}(1)$. Thus, if $m_n=n^\gamma$ for some $\gamma\in
(0,1)$ and
$r_n=n^\delta$ for some $\delta\in(\gamma/(s-1),\gamma)$, some
$s>2$, then \eqref{eqll} holds.
Condition (M1) is satisfied if $(1+s)^{-1}<\gamma<1/3$
and $\delta\in((1-\gamma)/s,\gamma)$. Thus \eqref{eqll} and (M1)
are always satisfied if $s$ can be chosen arbitrarily large.
\end{remark}
For our main result on the smoothed periodogram (see Theorem~\ref{thmsmoth}),
we finally need the condition:
\begin{longlist}[(M2)]
\item[(M2)]
The sequences $(m_n)$, $(r_n)$ from (M) also satisfy
the growth conditions
\[
 m_n^2 n \sum_{h=r_n+1}^{n}
\xi_h \to0,\qquad m_n r_n^3/n \to0.
\]
\end{longlist}
%
\begin{remark}\label{remrm2}
Condition (M2) is stronger than \eqref{eqll}. If $(X_t)$ is strongly
mixing with geometric or polynomial rate, a similar argument as
in Remark~\ref{remrm1} shows that (M2) holds for suitable choices of
$(r_n)$ and $(m_n)$.
\end{remark}
%
\subsection{The periodogram of extreme events}\label{subsecper}
In this section, we recall some of the results from Davis and
Mikosch \cite{davismikosch2009} concerning the estimation of the
spectral density $f_A$ defined in \eqref{eq5}.
Write
\begin{eqnarray*}
I_t=I_{\{X_t/a_m\in A\}} ,\qquad \wt I_t=
I_t-p_0 ,\qquad p_0=EI_t=P
\bigl(a_m^{-1} X\in A\bigr) ,\qquad t=1,\ldots,n
\end{eqnarray*}
for some sequence  $m=m_n\to\infty$ such that $m_n/n\to0$ as in
condition (M) above. We suppress the dependence of $I_t$ on $A$ and
$a_m$.
We introduce the estimators
%
\begin{eqnarray}
\label{eq16} I_{nA}(\la)=\frac{m_n} n \Biggl|\sum
_{t=1}^n \wt I_t \ex^{-\mathrm{i} t \la}
\Biggr|^2 ,\qquad \la\in[0,\uppi] \quad \mbox {and} \quad\wh P_m(A)=
\frac{m_n} n \sum_{t=1}^n
I_t.
\end{eqnarray}
It follows from Theorem 3.1 in
\cite{davismikosch2009}
that
%
\begin{equation}
\label{eq14} \wh P_m(A)=\frac{m_n} n \sum
_{t=1}^n I_t \stackrel{L^2}
{\to} \mu_0(A) =\lim_{\nto} m_n P
\bigl(a_m^{-1} X\in A\bigr) ,
\end{equation}
provided $A$ is a continuity set with respect to the limiting measure $\mu_0$.
The conditions $m_n\to\infty$ and $m_n/n\to0$ cannot be avoided
since we need that $E\wh P_m(A)=m_nP(a_m^{-1} X\in A)\to\mu_0(A)$ and
then we also get $\var(\wh P_m(A))=\mathrm{O}(m_n/n)$.

Davis and Mikosch \cite{davismikosch2009}, Theorem 5.1, also proved that
the \emph{lag-window estimator} or \emph{truncated periodogram}
%
\begin{equation}
\label{eqlagwindow} \wh f_{nA}(\la)=\wt\gamma_n(0)+2 \sum
_{h=1}^{r_n} \cos(\la h) \wt
\gamma_n(h)
\end{equation}
with $\wt\gamma_n(0)= (m/n) \sum_{t=1}^{n} I_t $ and $\wt\gamma
_n(h)= (m/n) \sum_{t=1}^{n-h} \wt I_t \wt I_{t+h}$, $h > 0$, for fixed
$\la\in(0,\uppi)$,
satisfies the relations
%
\begin{equation}
\label{eq15} E\wh f_{nA}(\la)\to\mu_0(A)
f_A(\la) \quad\mbox{and}\quad E \bigl(\wh f_{nA}(\la)-
\mu_0(A) f_A(\la) \bigr)^2\to0
\end{equation}
under condition (M), if $A$ is a $\mu_0$-continuity set and the sets
$A\times\ov\bbr_0^{k-1}\times A$ are continuity
sets with respect to $\mu_k$, $k\ge1$, and $m_n r_n^2=\mathrm{O}(n)$.
If we combine \eqref{eq14} and \eqref{eq15} we have
for fixed $\la\in(0,\uppi)$,
%
\begin{equation}
\label{eqstandardizedp} \frac{\wh f_{nA}(\la)}{\wh P_m(A)}\stp f_A(\la).
\end{equation}

A natural self-normalized estimator of the spectral density $
f_A(\la)$ in \eqref{eq5}
is
the following analog of the periodogram
\[
\wt I_{nA}(\la)= \frac{I_{nA}(\la)}{\wh P_m(A)} =\frac{ |\sum_{t=1}^n \wt I_t
\ex^{-\mathrm{i}t\la}  |^2}{\sum_{t=1}^n I_t} , \qquad\la\in[0,\uppi]
,
\]
In contrast to $\wh f_{nA}(\la)$ one does not need to know the
quantities $m_n$ and $r_n$ which appear
in the definition of $\wh f_{nA}(\la)$ and are hard to determine for
practical estimation purposes. We call $\wt I_{nA}(\la)$
the \emph{standardized periodogram}. However, we know from
theory for the classical periodogram of the stationary process
$(X_t)$,
given by
\[
J_{n,X}(\la)=n^{-1} \Biggl|\sum_{t=1}^n
X_t \ex^{-\mathrm{i}t\la}\Biggr |^2 ,\qquad \la\in[0,\uppi] ,
\]
that $J_{n,X}(\la)$ is
\emph{not} a consistent estimator of the spectral density $f_X(\la)$
of the
process
$(X_t)$ even in the case when $(X_t)$ is i.i.d. and has finite
variance; see, for example, Proposition 10.3.2 in Brockwell and Davis
\cite{brockwelldavis1991}. To achieve consistent estimation of
$f_X(\la)$ one needs
to truncate the periodogram, similarly to $\wh
f_{nA}(\la)$, or to apply smoothing techniques to neighboring
periodogram ordinates.
A similar observation applies to the periodogram for extremal events, $I_{n,A}(\la)$; see Section \ref{secprop}.

\section{Examples}\label{subsecexam}
In this section, we collect some examples of regularly varying stationary time series
models, give their extremograms \eqref{eq1} and, if possible,
give an explicit expression of the corresponding spectral density
\eqref{eq5}. However, in general,
the extremogram is too complicated
and one cannot calculate the
Fourier series \eqref{eq5}. Some of the examples below are taken
from Davis and Mikosch~\cite{davismikosch2009}.
\subsection{IID sequence}
Consider an i.i.d. real-valued sequence $(Z_t)$ such that
%
\begin{equation}
\label{eq6} P(Z > x) \sim p x^{-\alpha} L(x) \quad\mbox{and}\quad P(Z\le -x)\sim q
x^{-\alpha} L(x) ,\qquad \xto,
\end{equation}
where $\alpha>0$, $p,q\ge0$, $p+q=1$ and $L$ is a slowly varying function. It
is well
known (e.g., Resnick \cite{resnick1986,resnick1987}) that $(Z_t)$ is
regularly varying with index $\alpha$. The limiting measures $\mu_h$ are
concentrated on the axes:
\[
\mu_h(\mathrm{d}x_0,\ldots,\mathrm{d}x_h)=\sum
_{i=0}^h \la_\alpha(\mathrm{d}x_i)
\prod_{i\ne
j}\vep_0 \,{\mathrm{d}x_j}
,
\]
where $\vep_y$ denotes Dirac measure at $y$,
$\la_\alpha(x,\infty]=px^{-\alpha}$,
$\la_\alpha[-\infty,-x]=qx^{-\alpha}$, $x>0$.
Then for any $A$ bounded away from zero,
\[
\rho_{A}(h)=0 ,\qquad h\ge1 \quad \mbox{and}\quad f_A\equiv1.
\]
The conditions (M), (M1) and (M2) are trivially satisfied in this case.

\subsection{Stochastic volatility model} Let $(\sigma_t)$ be a
strictly stationary sequence of non-negative random variables with
$E\sigma^{\alpha+\delta}<\infty$ for some $\delta>0$, independent
of the i.i.d. regularly varying sequence $(Z_t)$ with index $\alpha>0$,
satisfying the tail balance
condition \eqref{eq6}. The process
\[
X_t=\sigma_t Z_t ,\qquad t\in\bbz,
\]
is a \emph{stochastic volatility process}. It is a regularly varying sequence with index $\alpha
$ and limiting
measures concentrated on the axes. The extremogram and the spectral
density coincide with these quantities in the i.i.d. case; see Davis
and Mikosch \cite{davismikosch2001}.
As discussed in Davis and Mikosch~\cite{davismikosch2009b},
the process $(X_t)$ inherits the strong mixing property and the same
rate function from the volatility process $(\sigma_t)$. In particular, if
$(\sigma_t)$ is strongly mixing with \emph{geometric rate}, $(X_t)$ is
also strongly mixing with
geometric rate, and then the conditions \eqref{eqll}, (M1) and (M2)
are satisfied; see Remarks~\ref{remrm1} and \ref{remrm2}.
Condition \eqref{eq33} also holds if $E\sigma^{4\alpha}<
\infty$; see Davis and Mikosch \cite{davismikosch2009}.

The situation of a vanishing $\rho_A$ is rather incomplete
information about tail dependence.
Hill~\cite{hill2011} proposed to use an alternative
lag-wise dependence measure of the form
$\lim_{\xto} P(X_h>x,X_0>x)/ [P(X_0>x)]^2 -1$ which in general does
not vanish. This measure is in agreement with the asymptotic tail independence
conditions of Ledford and Tawn~\cite{ledfordtawn2003}.

The mentioned literature \cite{davismikosch2001,davismikosch2009b}
focuses on stochastic volatility processes with i.i.d. regularly varying noise $(Z_t)$ with index
$\alpha$ and stochastic volatility satisfying the moment condition
$E\sigma^{\alpha+\delta}<\infty$ for some $\delta>0$. Mikosch and
Rezapur \cite{mikoschrezapur2012} consider regularly varying stochastic volatility processes with index $\alpha$ when
the sequence $(\sigma_t)$ is regularly varying with index $\alpha$,
$E|Z|^{\alpha+\delta}<\infty$ for some $\delta>0$ and they give examples
with $\rho_A\ne0$ and
$f_A\not\equiv1$ for $A$ bounded away from zero. The aforementioned
comments about mixing also apply in this setting.

\subsection{ARMA process}
Consider the linear process
%
\begin{equation}
\label{eq8} X_t = \sum_{j=0}^\infty
\psi_j Z_{t-j} , \qquad t\in\bbz,
\end{equation}
where $(Z_t)$ is an i.i.d. one-dimensional regularly varying sequence with index
$\alpha>0$ and tail
balance condition \eqref{eq6}. We choose the coefficients from the
ARMA equation $\psi(z)= 1+\sum_{i=1}^\infty\psi_i z^i=\theta
(z)/\phi(z)$, $z\in\bbc$, where
\begin{eqnarray*}
\phi(z)=1-\phi_1 z-\cdots-\phi_r z^r\quad
\mbox{and} \quad\theta(z)=1+\theta_1z +\cdots+ \theta_s
z^s
\end{eqnarray*}
for integers $r,s\ge0$,
and the coefficients $ \theta_i,\phi_i$ are chosen such that $\phi(z)$
and $\theta(z)$ have no common zeros and $\phi(z)\ne0$ for $|z|\le
1$.
It is well known that $X$ is regularly varying with index $\alpha$; see,
for example, Appendix A3.3 in Embrechts \textit{et al.}
\cite{embrechtskluppelbergmikosch1997} or Mikosch and Samorodnitsky
\cite{mikoschsamorodnitsky2000}. The proofs in the latter references
use the fact that $X_t^{(s)} = \sum_{j=0}^s \psi_j Z_{t-j}$, $s\ge1$,
is regularly varying as a simple consequence of the fact that linear
combinations of
i.i.d. regularly varying random variables are regularly varying; see
Feller \cite{feller1971},
page~278;
cf. Lemma 1.3.1 in
\cite{embrechtskluppelbergmikosch1997}. Moreover,
%
\begin{equation}
\label{eq9} \lim_{s\to\infty}\limsup_{\nto} n P
\bigl(a_n^{-1}\bigl|X_t-X_t^{(s)}\bigr|>
\vep \bigr)=0 ,\qquad \vep>0.
\end{equation}
Then it follows from Lemma 3.6 in Jessen and Mikosch
\cite{jessenmikosch2006} that $X_t$ is regularly varying.

The vector
$(X_0^{(s)}, \ldots,X_h^{(s)})$ is also regularly varying with
index $\alpha$. This fact follows from an application of a
multivariate version of
Breiman's lemma \cite{breiman1965} (see Basrak \textit{et~al.}~\cite{basrakdavismikosch2002})
or the fact that linear operations preserve regular variation; see
Lemma 4.6 in \cite{jessenmikosch2006}. Since
\eqref{eq9} holds a straightforward
multivariate extension of Lemma 3.6 in
\cite{jessenmikosch2006} yields that $(X_0,\ldots,X_h)$ is
regularly varying
for every $h\ge0$.

The same arguments leading to the asymptotic tail behavior of $X_t$
(see, e.g., Appendix A3.3 in Embrechts \textit{et al.}
\cite{embrechtskluppelbergmikosch1997},
Mikosch and Samorodnitsky \cite{mikoschsamorodnitsky2000}) yield for
$A=(1,\infty)$,
%
\begin{equation}
\label{eq10} \rho_{A}(h) = \frac{\sum_{i=0}^\infty [ p (\min(\psi_{i}^+ ,
\psi_{i+h}^+))^\alpha+ q (\min(\psi_{i}^- ,
\psi_{i+h}^-))^\alpha ]}{ \sum_{i=0}^{\infty}  [ p
(\psi_i^{+})^\alpha+ q (\psi_i^{-})^\alpha ]} ,\qquad h\ge1.
\end{equation}
This formula was given in \cite{davismikosch2009} for symmetric $Z$
when $p=q=0.5$.

Doukhan \cite{doukhan1994}, Theorem 6 on page~99,
shows that $(X_t)$ is $\beta$-mixing, hence strongly mixing,
with geometric rate if
$Z$ has a positive Lebesgue density in some neighborhood of the
expected value of $Z$ (provided it exists) and Pham and Tran
\cite{phamtran1985} proved the same statement under the condition that
$Z$ has a Lebesgue density and a finite $p$th moment for some $p>0$.
Hence \eqref{eqll}, (M1) and (M2) are
satisfied under these conditions; see Remarks~\ref{remrm1}
and \ref{remrm2}. Next, we verify condition \eqref{eq33}.
We observe that it trivially holds for an $s$-dependent sequence for any
integer \mbox{$s\ge1$.}
Hence, it is satisfied for any moving average of order $s$,
in particular for the truncated sequence $(X_t^{(s)})$.
For ease of presentation, we assume $\epsilon=1$.
Since $X_h^{(h-1)}$ and $X_0$ are independent we have
\begin{eqnarray*}
P\bigl(|X_h|>a_m\mid|X_0|>a_m\bigr)&
\le& P\bigl(\bigl|X_h^{(h-1)}\bigr|>0.5 a_m\bigr)+ P
\bigl(\bigl|X_h-X_h^{(h-1)}\bigr|>0.5 a_m,
\mid|X_{0}|>a_m\bigr)
\\
&\le&I_1+I_2.
\end{eqnarray*}
Recall that there exist $\varphi\in(0,1)$ such that $|\psi_i|\le
\varphi^i$
for $i$
sufficiently large; see Brockwell and Davis \cite
{brockwelldavis1991}, Chapter 3.
We have for a positive constant $c>0$, for every $k\ge1$,
\begin{eqnarray*}
\sum_{h=k+1}^{r_n} I_1&\le&
r_n P \Biggl(\sum_{i=0}^\infty|
\psi_{i}|| Z_i|>0.5 a_m \Biggr)
\\
&\sim& c r_n P\bigl(|Z|>a_m\bigr)= \mathrm{o}(1) \qquad\mbox{as $\nto$.}
\end{eqnarray*}
(Here and in what follows, $c$ denotes any constant whose value is not
of interest.)
For sufficiently large $k$, we have in view of the uniform convergence theorem for regularly varying functions (see Bingham \textit{et al.} \cite
{binghamgoldieteugels1987},
Section 1.2),
\begin{eqnarray*}
\sum_{h=k+1}^{r_n} I_2&\le&c
m_n \sum_{h=k+1}^{r_n} P \Biggl(
\sum_{i=h+1}^\infty|\psi_{i}||
Z_i|>0.5 a_m \Biggr)
\\
&\le& c m_n \sum_{h=k+1}^{r_n} P
\Biggl(\varphi^h\sum_{i=0}^\infty
\varphi^i | Z_i|>0.5 a_m \Biggr)
\\
&\le&c\sum_{h=k+1}^{r_n}
\varphi^{\alpha h}\le c \varphi^{\alpha
(k+1)}/\bigl(1-\varphi^\alpha
\bigr) ,
\end{eqnarray*}
and the right-hand side converges to zero as $\kto$. Thus we proved that
(M), (M1) and (M2) hold for ARMA processes if the noise has some
Lebesgue density.

If $\var(X)<\infty$ relation \eqref{eq10} bears some similarity
with the autocorrelation function of
$(X_t)$ given by $\rho(h)=\sum_{i=1}^\infty
\psi_i\psi_{i+h}/\sum_{i=1}^\infty\psi_i^2$. Replacing $\rho_{A}$ in
\eqref{eq5} by $\rho$, one obtains the well-known spectral density of
a causal ARMA process (up to a constant multiple):
$f_X(\la)=(2\uppi)^{-1}|\theta(\ex^{-\mathrm{i}\la})|^2/|\phi(\ex^{-\mathrm{i}\la})|^2$,
$\la\in[0,\uppi]$. Such a compact formula can in general not be
derived for $ f_A$.
An exception is a causal $\operatorname{ARMA}(1,1)$ process; see Section~\ref{exam1}.
There are various analogies between the functions $\rho$ and
$\rho_A$ for causal invertible ARMA processes. In this case, $\psi
_h\to0$ as $h\to\infty$ at an exponential rate and therefore both
$\rho(h)$ and $\rho_A(h)$
decay exponentially fast to zero as well. The latter property also
makes the
spectral densities $f_X$ and $f_A$ analytical functions bounded
away from infinity. We also mention that for an $\operatorname{MA}(q)$ process,
$\rho(h)=\rho_A(h)=0$ for $h>q$.

\subsection{Max-moving averages} Consider a regularly varying i.i.d.
sequence $(Z_t)$ with index $\alpha>0$ and tail balance parameters~$p,q$; see
\eqref{eq6}. For a real-valued sequence $(\psi_j)$, the process
%
\begin{equation}
\label{eq13} X_t = \bigvee_{i=0}^{\infty}
\psi_i Z_{t-i},\qquad t\in\bbz,
\end{equation}
is a \emph{max-moving average}. We will also assume that $|\psi_j|\le
c$, $j\ge0$, for some constant $c$ and $\psi_0=1$.
Obviously, if $X$ is finite a.s.,
$(X_t)$ constitutes a strictly stationary process. The random variable $X$ does
not assume the value $\infty$ if $\lim_{\xto}P(X> x)=0$.
We have
\[
P(X>x)= P \Biggl(\bigvee_{i=0}^{\infty}
\psi_i Z_i>x \Biggr)= 1-\lim_{\nto}
\prod_{i=0}^n P(\psi_i Z\le x)
.
\]
The product $\prod_{i=0}^\infty P(\psi_i Z\le x)$ converges
if $\sum_{i=0}^\infty P(\psi_i Z >x)<\infty$.
By regular variation of $Z$,
this amounts to the condition
\[
\psi_+=\sum_{i=0}^\infty\bigl[p \bigl(
\psi_i^+\bigr)^\alpha+ q \bigl(\psi_i^-
\bigr)^\alpha\bigr]<\infty.
\]
A Taylor expansion and regular variation of $Z$ yield
%
\begin{equation}
\label{eq11} P(X>x)= 1- \ex^{ -(1+\mathrm{o}(1))P(|Z|>x)\psi_+}\sim P\bigl(|Z|>x\bigr)\psi_+\to0 ,\qquad \xto.
\end{equation}
We also have $P(X\le-x)=\mathrm{O}(P(|Z|>x))$.
Hence, $X$ is regularly varying with index $\alpha$
if $0<\psi_+<\infty$. We always assume the latter condition.

We show that $(X_t)$ is regularly varying. Consider the
truncated max-moving average process for $s\ge0$,
\[
X_t^{(s)} = \bigvee_{i=0}^s
\psi_i Z_{t-i} ,\qquad t \in\bbz.
\]
Regular variation of $(X_0^{(s)},\ldots,X_h^{(s)})$ is a consequence of
regular variation of $(Z_t)$ and the fact that regular variation is
preserved under
the max-operation acting on independent components. Moreover,
\begin{eqnarray*}
\lim_{s\to\infty}\limsup_{\nto} n P
\Biggl(a_n^{-1} \bigvee_{i=s+1}^\infty
\psi_i Z_{t-i}>x \Biggr)&=&c \lim_{s\to\infty}
\sum_{i=s+1}^\infty \bigl[p \bigl(
\psi_i^+\bigr)^\alpha+ q \bigl(\psi_i^-
\bigr)^\alpha\bigr]=0.
\end{eqnarray*}
Then an application of Lemma 3.6 in Jessen and Mikosch
\cite{jessenmikosch2006} shows that $(X_0,\ldots,X_h)$ is regularly varying
with index $\alpha$ for every $h\ge0$.

Next, we determine the extremogram $\rho_{A}$ corresponding to the
set $A=(1,\infty)$. For $h\ge1$, we have
\begin{eqnarray*}
P(X_h>x,X_0>x)&=&P \Biggl(\bigvee
_{i=0}^\infty\psi_i Z_{-i}>x ,
\bigvee_{i=-h}^{-1} \psi_{i+h}Z_{-i}
\vee\bigvee_{i=0}^\infty\psi _{i+h}
Z_{-i}>x \Biggr)
\\
&=&P \Biggl(\bigvee_{i=0}^\infty(
\psi_i Z_{-i}) \wedge(\psi _{i+h}Z_{-i})
>x \Biggr) +\mathrm{o}\bigl(P\bigl(|Z|>x\bigr)\bigr)
\\
&\sim&P\bigl(|Z|>x\bigr) \sum_{i=0}^\infty \bigl[ p
\bigl(\min\bigl(\psi_{i}^+ , \psi _{i+h}^+\bigr)
\bigr)^\alpha + q \bigl(\min\bigl(\psi_{i}^- ,
\psi_{i+h}^-\bigr)\bigr)^\alpha \bigr].
\end{eqnarray*}
Finally, in view of \eqref{eq11}, $\rho_{A}(h)$ is given by
\eqref{eq10}, that is, the linear
process \eqref{eq8} and the max-moving average \eqref{eq13} have the
same extremogram provided the coefficients $(\psi_j)$ and
the distribution of $Z$ are the same. Hence, their spectral densities
$ f_A$ are the same as well.\looseness=-1

As for ARMA processes, mixing conditions for infinite
max-moving processes are not easily verified and
additional conditions on the noise $(Z_t)$ are needed. Assume that
$(Z_t)$ is i.i.d. with common
Fr\'echet distribution $\Psi_\alpha(x)=\ex^{-x^{-\alpha}}$, $x>0$,
for some
$\alpha>0$. Then
$(X_t)$ constitutes a stationary max-stable process. For such processes,
Dombry and Eyi-Minko \cite{dombryeyi-minko2012}
proved rather general sufficient conditions for $\beta$-mixing,
implying strong mixing. An application of their Corollary 2.2 implies
that the condition $|\psi_h|\le c_0 \ex^{-c_1 h}$, $h\ge1$,
for suitable constants $c_1,c_2>0$ implies strong mixing of $(X_t)$
with geometric rate function $(\xi_h)$.
In this situation, (M), (M1) and~(M2) are satisfied.

\section{Basic properties of the periodogram}\label{secprop}
In this section, we study some basic properties of the periodogram $I_{nA}(\la)$ for extremal events defined in \eqref{eq16}. Notice that
\[
I_{nA}(\la)=\tfrac12 \bigl[\bigl(\alpha_n(\la)
\bigr)^2+\bigl(\beta_n(\la)\bigr)^2 \bigr],
\]
where $\alpha_n(\la)$ and $\beta_n(\la)$ denote the normalized and
centered cosine
and sine transforms of $(I_t)_{t=1,\ldots,n}$:
\begin{eqnarray*}
\alpha_n (\lambda)&=& \biggl(\frac{2m_n}n
\biggr)^{1/2} \sum_{t=1}^{n} \wt
I_t \cos(\lambda t) ,
\\
\beta_n (\lambda)&=& \biggl( \frac{2m_n} n
\biggr)^{1/2} \sum_{t=1}^{n} \wt
I_t \sin(\lambda t).
\end{eqnarray*}
Here we suppress the dependence of $\alpha_n$ and $\beta_n$ on $a_m$
and the
set $A$ which is bounded away from zero.
For practical purposes, the periodogram will typically be evaluated at some
Fourier frequencies $\la= 2\uppi j/n$ for some
integer $j$.
If $\la\in(0,\uppi)$ is such a \emph{Fourier frequency,} then
\[
\sum_{t=1}^{n} \ex^{\mathrm{i} \la t}=0 ,
\]
and therefore the
$I_t$'s in $\alpha_n (\lambda)$ and $\beta_n (\lambda)$ are
automatically centered by their (in general unknown)
expectations $EI_t=p_0=P(a_m^{-1} X\in
A)$.

\subsection{The periodogram ordinates at
distinct frequencies are asymptotically uncorrelated}\label{subsecperiod}
Our first result is an analog of the fact that the sine and cosine
transforms of a stationary sequence at distinct fixed or Fourier
frequencies in
$(0,\uppi)$ are asymptotically uncorrelated.

\begin{proposition}\label{prop1}
Consider a strictly stationary $\bbr^d$-valued sequence $(X_t)$ which is regularly varying with index $\alpha>0$ and
satisfies the
mixing condition \textup{(M)}. Let $A\subset\ov\bbr^d_0$ be bounded
away from zero such that $A$ is a continuity set with respect to $\mu
_0$ and
$A \times\ov\bbr_0^{dh}$ and $A\times\ov\bbr_0^{d(h-1)} \times
A$ are continuity sets with respect to the limiting measures $\mu
_{h}$ for every
$h\ge1$; see Section~\ref{subsecregvar}. Also assume that
$\sum_{h\ge1}\rho_A(h)<\infty$.
Let $\lambda,\omega$ be either any two Fourier
or fixed frequencies
in $(0,\uppi)$.
\begin{longlist}[(3)]
\item[(1)]
If $\la,\w$ are distinct then
the covariances of the pairs
$(\alpha_n(\lambda),\beta_n(\omega))$,
$(\alpha_n(\lambda),\alpha_n(\omega))$,
$(\beta_n(\lambda),\beta_n(\omega))$ converge to zero as $n\to
\infty$.
\item[(2)]
The covariance of $(\alpha_n(\lambda),\beta_n(\la))$
converges to zero as $\nto$.
\item[(3)]
If $\la\in(0,\uppi)$ is fixed and if $(\la_n)$ are Fourier frequenc
ies
such that $\la_n\to\la,$
then the asymptotic variances are given by
\begin{eqnarray*}
\var\bigl(\alpha_n(\la_n)\bigr)&\sim&\var\bigl(
\alpha_n(\la)\bigr)\sim \var\bigl(\beta_n(
\la_n)\bigr)\sim\var\bigl(\beta_n(\la)\bigr)
\\
&\sim&\mu_0(A) \Biggl[1+2\sum_{h=1}^\infty
\cos(\la h)\rho_{A}(h) \Biggr]=\mu_0(A)
f_A(\la).
\end{eqnarray*}
\end{longlist}
\end{proposition}
%
\begin{remark}
The smoothness condition on the set $A$ ensures that the extremogram
$\rho_A$ with respect to $A$ is well defined; see Section~\ref
{subsecregvar}.
\end{remark}
%
\begin{remark}\label{rem2}
Since $E\alpha_n(\la)=E\beta_n(\la)=0$
an immediate consequence of part (3) is that
\begin{eqnarray*}
EI_{nA}(\la)= \frac12 \bigl[\var\bigl(\alpha_n(\la)
\bigr)+\var\bigl(\beta_n(\la )\bigr) \bigr] \sim\mu_0(A)
\Biggl[1+2\sum_{h=1}^\infty\cos(\la h)
\rho_{A}(h) \Biggr]=\mu_0(A)f_A(\la).
\end{eqnarray*}
Following the lines of the proof below, one can see that
the error one encounters in the above approximation is uniform for
$\la\in[a,b]\subset(0,\uppi).$ The same remark applies to the
quantities $EI_{nA}(\la_n)$ evaluated at
Fourier frequencies $\la_n\to\la\in(0,\uppi)$.
\end{remark}
\begin{pf} We start by calculating the asymptotic covariances.
Any of the
covariances can be written in the form
\begin{eqnarray*}
J&=& \frac{2m_n}{n}E \Biggl[ \sum_{s=1}^n
\sum_{t=1}^n \bigl(I_s
I_t-p_0^2\bigr) f_1 (
\lambda s)f_2 (\omega t) \Biggr]
\\
&=&\frac{2m_n}{n} \biggl(\sum_{1\le t=s\le n}+ \sum
_{1\le s\neq t\le n} \biggr) \bigl(p_{|s-t|}-p_0^2
\bigr) f_1 (\lambda s) f_2 (\omega t)
\\
&=&J_1+J_2 ,
\end{eqnarray*}
where $f_1$ and $f_2$ are cosine or sine functions
and
\[
p_{|t-s|}=P\bigl(a_m^{-1} X_s\in A,
a_m^{-1} X_t\in A\bigr) \qquad\mbox{for any
$s,t$.}
\]
We estimate $J_1$ separately for each possible combination of sine and cosine
functions $f_1,f_2$. We start with
$f_1(x)=\cos x$ and $f_2(x)=\sin x$.
Then, if $\la,\w$ are Fourier frequencies, so are
$\la\pm\w$ and therefore
\begin{eqnarray*}
J_1 &=&\bigl(p_0-p_0^2
\bigr) \frac{2m_n}{n} \sum_{t=1}^n
\cos(\lambda t)\sin (\omega t)
\\
& =&\bigl(p_0-p_0^2\bigr)
\frac{m_n}{n}\sum_{t=1}^n \bigl[ \sin
\bigl((\lambda +\omega) t\bigr)- \sin\bigl((\omega-\lambda)t\bigr) \bigr]=0.
\end{eqnarray*}
If $\la,\w$ are fixed frequencies, we conclude from \eqref{eqb}
that the sum on the right-hand side is bounded. Hence, $J_1=\mathrm{O}(n^{-1})$.

For $f_1(x)=f_2(x)=\cos x,$ we get
\begin{eqnarray*}
J_1 &=&\bigl(p_0-p_0^2
\bigr)\frac{2m_n}{n} \sum_{t=1}^n
\cos(\lambda t)\cos (\omega t)
\\
& =& \bigl(p_0-p_0^2\bigr)
\frac{m_n}{n} \sum_{t=1}^n \bigl[
\cos \bigl((\lambda+\omega) t\bigr)+ \cos\bigl((\omega-\lambda)t\bigr) \bigr].
\end{eqnarray*}
If $\la,\w$ are Fourier frequencies, so are $\la\pm\w$ and then the
right-hand side vanishes
unless $\la+\w= \uppi$. However, if
$\la+\w= \uppi$ the second sum vanishes and the first sum
is bounded. Therefore, $J_1=\mathrm{O}(n^{-1})$.
If $\la\ne\w$ are fixed it follows from \eqref{eqa} that the sum on
the right-hand side is bounded and therefore $J_1=\mathrm{O}(n^{-1})$.

For $f_1(x)=f_2(x)=\sin x$ we have
\begin{eqnarray*}
J_1 &=& \bigl(p_0-p_0^2
\bigr) \frac{2m_n}{n} \sum_{t=1}^n
\sin(\lambda t) \sin(\omega t)
\\
&=& \bigl(p_0-p_0^2\bigr)
\frac{m_n}{n} \sum_{t=1}^n \bigl[
\cos\bigl((\lambda -\omega) t\bigr)- \cos\bigl((\la+\w)t\bigr) \bigr].
\end{eqnarray*}
The same arguments as above show that $J_1=\mathrm{O}(n^{-1})$
both for Fourier and fixed frequencies $\la\ne\w$.

Next, we consider $J_2$. We start with
$\cov(\alpha_n(\la),\beta_n(\la))$. If $\la$ is a Fourier frequency,
we have $\sin(\la n)=0$. Hence, by \eqref{eqe},
\begin{eqnarray*}
J_2&=& \frac{2m_n}{n} \sum_{h=1}^{n-1}
\bigl(p_{h}-p_0^2\bigr)\sum
_{s=1}^{n-h} \bigl[\sin(\lambda s) \cos\bigl(\la(s+h)
\bigr)+ \cos(\lambda s) \sin\bigl(\la( s+h)\bigr) \bigr]
\\
&=&- \frac{ 2m_n}{n} \sum_{h=1}^{n-1}
\bigl(p_{h}-p_0^2\bigr)\sin(\la h).
\end{eqnarray*}
By definition of strong mixing,
$|p_h-p_0^2|\le\xi_h$. Then, by condition (M),
\[
|J_2|\le\frac{2m_n}{n} \sum_{h=1}^\infty
\xi_h=\mathrm{O}(m_n/n).
\]
The same argument applies for a fixed frequency $\la$ since
the expressions in \eqref{eqe} are bounded for every $n$ and $h<n$.

If $\la\ne\w$ are fixed frequencies, we conclude from
\eqref{eqf}--\eqref{eqh} and condition (M)
that there exist constants
$c(\lambda,\omega)$ such that
\begin{eqnarray*}
|J_2| & =& \Biggl| \frac{2m_n}{n} \sum_{h=1}^{n-1}
\bigl(p_{h}-p_0^2\bigr)\sum
_{s=1}^{n-h} \bigl(f_1(\lambda s)
f_2\bigl(\omega (s+h)\bigr)+ f_1\bigl(\lambda(s+h)
\bigr)f_2(\omega s) \bigr)\Biggr |
\\
&\le& c(\lambda,\omega) \frac{m_n}{n} \sum_{h=1}^{n-1}
\bigl|p_{h}-p_0^2\bigr| \leq c(\lambda,\omega)
\frac{m_n}{n} \sum_{h=1}^{\infty}\xi
_h=\mathrm{O}(m_n/n).
\end{eqnarray*}
Now we consider the case of two distinct Fourier frequencies $\la,\w$.
We start with $f_1(x)=\cos x$ and $f_2(x)=\sin x$.
If $\la+\w$ and $|\la-\w|$ are bounded away from zero, we can use the
argument for general distinct frequencies. Now assume that $\la+\w\le
0.1$ say.
Since $\la,\w$ are Fourier frequencies a glance
at \eqref{eqf}--\eqref{eqh} shows that one has to find suitable
bounds for
\[
\frac{|\sin((n-h+1)(\la+\w)/2)|}{|\sin((\la+\w)/2)|} = \frac{|\sin((-h+1)(\la+\w)/2)|}{|\sin((\la+\w)/2)|}.
\]
If $h (\la+\w)\le0.1$ Taylor
expansions for the nominator and the denominator show that
the right-hand side is bounded by $c h$. If $h (\la+\w)>0.1$ bound the
nominator by
1 and Taylor expand the denominator to conclude that the right-hand
side is
bounded by
$c h$ for some constant $c>0$ as well.
Then, by \eqref{eqf}, for fixed $k$,
\[
|J_2|  \le c \Biggl[\frac{m_n}{n} \sum
_{h=1}^{k}\bigl|p_h-p_0^2\bigr|
+ m_n \sum_{h=k+1}^{r_n}\bigl|p_h-p_0^2\bigr|+
m_n \sum_{h=r_n+1}^\infty\xi
_h \Biggr].
\]
The right-hand side vanishes by virtue of condition (M),
first letting $\nto$ and then $\kto$. The case of small $|\la-\w|$,
$|\la-\w|\le0.1$
say, can be treated analogously.

The remaining cases $f_1(x)=f_2(x)=\cos x$ and $f_1(x)=f_2(x)=\sin x$
can be treated in the same way by exploiting \eqref{eqg} and
\eqref{eqh}.

Now we turn to the asymptotic variances. We restrict ourselves
to $\alpha_n(\la)$ for fixed $\la\in(0,\uppi)$; the variance of
$\beta_n(\la)$ and the case of Fourier frequencies can be treated
analogously. Write

We have
\begin{eqnarray*}
\var\bigl(\alpha_n(\la)\bigr)&=&\frac{2m_n}{n} \Biggl[
\bigl(p_0-p_0^2\bigr) \sum
_{t=1}^n \bigl(\cos(\la t)\bigr)^2 + 2
\sum_{h=1}^{n-1}\bigl(p_h-p_0^2
\bigr)\sum_{t=1}^{n-h} \cos(\la t)\cos\bigl(
\la(t+h)\bigr) \Biggr].
\end{eqnarray*}
For any frequency $\la\in(0,\uppi)$ bounded away from zero and $\uppi$, the
relation\break  $n^{-1}\sum_{t=1}^n (\cos(\la t))^2\sim
0.5$ holds. Moreover, $\cos(\la t)\cos(\la(t+h))=0.5[\cos(\la
h) +\break \cos(\la(2t+h))]$. Similar calculations as above yield
\begin{eqnarray*}
\var\bigl(\alpha_n(\la)\bigr)&\sim& m_n p_0
+ 2 m_n \sum_{h=1}^{n-1}
\bigl(p_h-p_0^2\bigr) (1-h/n) \cos(\la h)
\\
&\sim&\mu_0(A) \Biggl[1+ 2 \sum_{h=1}^{\infty}
\rho_{A}(h)\cos(\la h) \Biggr].
\end{eqnarray*}
This concludes the proof.
\end{pf}
%
\subsection{Central limit theorem}\label{subsecclt}
Our next result shows that the periodogram ordinates at distinct
frequencies are asymptotically independent and exponentially distributed.
%
\begin{theorem}\label{thm1} Consider a strictly stationary $\bbr^d$-valued
sequence $(X_t)$ which
is regularly varying with index $\alpha>0$. Let $A \subset\ov\bbr
_0^d$ be
bounded away satisfying the smoothness conditions of Proposition~\ref{prop1}.
Assume that the conditions \emph{(M)}, \emph{(M1)} and
$\sum_{h\ge1}\rho_A(h)<\infty$
hold. Consider any fixed frequencies
$0<\lambda_1<\cdots< \lambda_N <\uppi$ for some $N\ge1$.
Then the following central limit theorem  holds:
%
\begin{equation}
\label{eq90} \bfZ_n= \bigl(\alpha_n(\la_i),
\beta_n(\la_i) \bigr)_{i=1,\ldots
,N}\std \bigl(\alpha(
\la_i),\beta(\la_i) \bigr)_{i=1,\ldots,N} ,\qquad \nto,
\end{equation}
where the limiting vector has $N({\mathbf{0}},\Sigma_N)$ distribution with
\[
\Sigma_N=\mu_0(A) \operatorname{diag} \bigl(f_A(
\la_1),f_A(\la_1),\ldots
,f_A(\la_N),f_A(\la_N)
\bigr).
\]
The limit relation \eqref{eq90} remains valid if the frequencies
$\la_i$, $i=1,\ldots,N$, are replaced by distinct Fourier frequencies
$\w_i(n)\to\la_i\in(0,\uppi)$ as
$\nto$. The limits $\la_i$ do not have to be distinct.
\end{theorem}
Then the following result is immediate.
%
\begin{corollary}\label{cor1}
Assume the conditions of Theorem~\ref{thm1}.
Let $(E_i)$ be a sequence of i.i.d. standard exponential random variables.
\begin{enumerate}
\item
Consider any fixed frequencies
$0<\lambda_1<\cdots< \lambda_N <\uppi$ for some $N\ge1$.
Then the following relations hold:
\begin{eqnarray*}
\bigl(I_{nA}(\lambda_i) \bigr)_{i=1,\ldots,N}&\std&
\mu_0(A) \bigl(f_A(\la_i) E_i
\bigr)_{i=1,\ldots,N} ,\qquad \nto ,
\\
\bigl(\wt I_{nA}(\lambda_i) \bigr)_{i=1,\ldots,N}&\std&
\bigl(f_A(\la_i) E_i \bigr)_{i=1,\ldots,N}
, \qquad \nto.
\end{eqnarray*}
\item
Consider any distinct Fourier frequencies $\w_i(n)\to\la_i\in(0,\uppi)$
as $\nto$, $i=1,\ldots,N$. The limits $\la_i$ do not have to be
distinct. Then
the following relations hold:
\begin{eqnarray*}
\bigl(I_{nA}\bigl(\w_i(n)\bigr) \bigr)_{i=1,\ldots,N}&
\std& \mu_0(A) \bigl(f_A(\la_i)
E_i \bigr)_{i=1,\ldots,N} , \qquad\nto ,
\\
\bigl(\wt I_{nA}\bigl(\w_i(n)\bigr) \bigr)_{i=1,\ldots,N}&
\std& \bigl(f_A(\la_i) E_i
\bigr)_{i=1,\ldots,N} ,\qquad  \nto.
\end{eqnarray*}
\end{enumerate}
\end{corollary}
\begin{pf*}{Proof of the Theorem~\ref{thm1}} We will prove \eqref{eq90}
by applying
the Cram\'er--Wold device, that is, we will show that for any choice of
constants $\bfc\in\bbr^{2N}$,
%
\begin{equation}
\label{eqclt} \bfc'\bfZ_n \std N\bigl(0,
\bfc'\Sigma_N \bfc\bigr).
\end{equation}
The proof of the result for distinct converging Fourier frequencies is analogous
and therefore omitted.
We will prove \eqref{eqclt} by applying
the method of small and large blocks. The difficulty we encounter
here is that, due to the presence of sine and cosine functions, we
are dealing with partial sums of non-stationary sequences.
For $t=1,\ldots,n$, we write
%
\begin{equation}
\label{eq13a} Y_{nt}= \biggl(\frac{2m_n}n \biggr)^{1/2}
\wt I_t \sum_{j=1}^{N}
\bigl[c_{2j-1}\cos(\lambda_j t)+ c_{2j}\sin(
\la_j t)\bigr] ,\qquad t=1,\ldots ,n.
\end{equation}
For ease of presentation, we always assume that $n/m_n =k_n$ is an integer;
the general case can be treated in a similar way.
Consider the large blocks
\[
K_{ni}=\bigl\{(i-1)m_n+1,\ldots,im_n\bigr\},\qquad
i=1,\ldots,k_n ,
\]
the index sets
$\wt K_{ni}$, which consist of all but the first $r_n$ elements of
$K_{ni}$, and the small blocks
$J_{ni}=K_{ni}\setminus\wt K_{ni}$.
In view of condition (M), $r_n/m_n \to0$ and $m_n\to\infty$, the
sets $\wt K_{ni}$ and
$J_{ni}$ are non-empty for large $n$.
For any set $B\subset\{1,\ldots,n\}$, we write
$
S_n(B)=\sum_{t\in B} Y_{nt}
$.
First, we show that the joint contribution of the sums over the small blocks
to $\bfc'\bfZ_n$ is asymptotically negligible.
\end{pf*}
%
\begin{lemma}Under the conditions of Theorem~\ref{thm1}, the following
relation holds:
%
\begin{equation}
\label{eqA11} \var \Biggl(\sum_{i=1}^{k_n}
S_n (J_{ni}) \Biggr)\to0, \qquad\nto.
\end{equation}
\end{lemma}
\begin{pf}
We have
\begin{eqnarray*}
\var \Biggl(\sum_{i=1}^{k_n}
S_n (J_{ni}) \Biggr) &\leq& \sum
_{i=1}^{k_n} \var\bigl(S_n(J_{ni})
\bigr)+2\sum_{1\le i_1<i_2\le k_n} \bigl|\cov \bigl(S_n(J_{ni_1}),S_n(J_{ni_2})
\bigr)\bigr|
\\
&=&P_1+P_2.
\end{eqnarray*}
Due to the sum structure of $Y_{nt}$ given in \eqref{eq13a}
each of the sums $S_{n}(J_{ni})$ can be written as a sum of $2N$
subsums where each of these subsums only involves either the functions
$\cos(\la_jt)$ or $\sin(\la_j s)$ for some $j\le N$.
Then each of the terms
$\var(S_n(J_{ni}))$ and $|\cov (S_n(J_{ni_1}),S_n(J_{ni_2}) )|$
is bounded by a linear combination of the variances/covariances
of such subsums. In other words, it suffices to prove \eqref{eqA11}
for $N=1$. We give the corresponding calculations only for the
functions $\cos(\la t)$ where $\la$ stands for any of the frequencies
$\la_j$.
The calculations are similar to those in the proof of Proposition~\ref{prop1}.
For any $i\le k_n$ and fixed $k\ge1$, condition (M) ensures that
there is a constant $c(k)$ such that for large $n$,
\begin{eqnarray*}
\var\bigl(S_n (J_{ni})\bigr) & =& \frac{2m_n}{n}
\Biggl[ \sum_{s=(i-1)m_n+1}^{(i-1)m_n +r_n}
\var(I_s) \bigl(\cos(\lambda s)\bigr)^2
\\
&&\hspace*{25pt}{}+2\sum_{h=1}^{r_n -1} \sum
_{s=(i-1)m_n+1}^{(i-1)m_n +r_n-h}\cov ( I_s,I_{s+h})
\cos(s\lambda) \cos\bigl(\lambda(s+h)\bigr) \Biggr]
\\
& \leq&\frac{2m_n}{n} \Biggl( r_n\bigl(p_0-p_0^2
\bigr) + 2 \sum_{h=1}^{r_n -1} (r_n
-h)\bigl|p_h-p_0^2\bigr| \Biggr)
\\
& \leq& c\frac{r_n}{n} \Biggl( m_n \sum
_{h=0}^{k}p_h + m_n \sum
_{h=k+1}^{r_n} p_h \Biggr)\le c(k)
(r_n/n) ,
\end{eqnarray*}
and the right-hand side does not depend on $i$.
Consequently, $P_1\le c(k) k_n r_n/n=c(k) r_n/m_n\to0$ for
every fixed $k$.
Similarly, for $i_1< i_2$,
\begin{eqnarray*}
&&\bigl|\cov\bigl(S_n (J_{ni_1}), S_n
(J_{n,i_2})\bigr)\bigr|\\
&&\quad = \frac{2m_n}{n} \Biggl| \Biggl[ \sum
_{s=(i_1-1)m_n+1}^{(i_1-1)m_n+r_n} \sum_{t=(i_2-1)m_n+1}^{(i_2-1)m_n+r_n}
\cov(I_t,I_s) \cos(\lambda s) \cos(\lambda t) \Biggr] \Biggr|
\\
& &\quad\leq c \frac{m_n}{n} \sum_{q=(i_2-i_1) m_n -(r_n-1)}^{(i_2-i_1)m_n
+r_n-1}
\bigl(r_n -\bigl|q- (i_2-i_1)m_n\bigr|
\bigr) \bigl|p_q -p_0^2\bigr|
\\
&&\quad \leq c\frac{ m_n r_n}{n} \sum_{q=(i_2-i_1)m_n
-(r_n-1)}^{(i_2-i_1)m_n +r_n-1}
\xi_q ,
\end{eqnarray*}
where $(\xi_t)$ is the mixing rate function.
Hence for large $n$, in view of condition (M),
\begin{eqnarray*}
|P_2| & \leq& c \frac{m_n r_n}{n}\sum_{i_1=1}^{k_n}
\sum_{i_2=i_1+1}^{k_n} \sum
_{q=(i_2 -i_1)m_n -(r_n-1)}^{(i_2 -i_1)m_n
+r_n-1} \xi_q
\\
& \leq& c \frac{m_n r_n}{n}\sum_{i_1=1}^{k_n -1}
\sum_{q=m_n+1-r_n}^\infty\xi_q \leq c
r_n \sum_{q=r_n+1}^ \infty
\xi_q=\mathrm{o}(1).
\end{eqnarray*}
This proves \eqref{eqA11}.
\end{pf}

Relation \eqref{eqA11} implies that $\bfc'\bfZ_n$ and
$\sum_{i=1}^{k_n} S_n (\wt K_{ni})$ have the same limit distribution
provided such a limit exists. Let $\wt S_n(\wt K_{ni})\eqd S_n(\wt
K_{ni})$ for $i=1,\ldots,k_n$ and assume
that $(\wt S_n(\wt K_{ni}))_{i=1,\ldots,k_n}$ has independent components.
A telescoping sum argument yields
\begin{eqnarray*}
&&\Biggl| E \prod_{l=1}^{k_n}
\ex^{\mathrm{i}t S_n(\wt K_{nl})} - E \prod_{s=1}^{k_n}
\ex^{\mathrm{i}t \wt S_n (\wt K_{ns})} \Biggr|
\\
&&\quad =  \Biggl| \sum_{l=1}^{k_n} E \Biggl[ \bigl(
\ex^{\mathrm{i}t S_n(\wt K_{nl})}- \ex^{\mathrm{i}t \wt S_n (\wt K_{nl})} \bigr) \prod_{s=1}^{l-1}
\ex^{\mathrm{i}t \wt S_n (\wt K_{ns})} \prod_{s=l+1}^{k_n}
\ex^{\mathrm{i}t S_n(\wt K_{ns})} \Biggr] \Biggr|
\\
&&\quad \leq \sum_{l=1}^{k_n} \Biggl| E \Biggl(\prod
_{s=1}^{l-1} \ex^{\mathrm{i}t \wt S_n (\wt K_{ns})} \bigl(
\ex^{\mathrm{i}t S_n(\wt K_{nl})}- \ex^{\mathrm{i}t \wt S_n (\wt K_{nl})} \bigr) \prod_{s=l+1}^{k_n}
\ex^{\mathrm{i}t S_n(\wt K_{ns})} \Biggr) \Biggr|
\\
&&\quad \leq4 k_n \xi_{r_n} \rightarrow0.
\end{eqnarray*}
In the last step, we used Theorem~17.2.1 in Ibragimov and Linnik
\cite{ibragimovlinnik1971} and condition (M1).
Hence, $\sum_{l=1}^{k_n} S_n (\wt K_{nl})$ and $\sum_{l=1}^{k_n} \wt
S_n (\wt K_{nl})$ have the same limits in distribution provided these
limits exist. In view of \eqref{eqA11} and the last conclusion the
central limit theorem \eqref{eqclt}
holds if and only if the same limit relation holds for $\sum_{i=1}^{k_n} \wt
S_n (K_{ni})$,
where $\wt S_n(K_{ni})\eqd S_n(K_{ni})$ and $(\wt S_n
(K_{ni}))_{i=1,\ldots,k_n}$ has independent components. Thus, we may
apply a classical central limit theorem for triangular arrays of independent random variables; see, for example,
Theorem 4.1 in Petrov \cite{petrov1995}.

According to this result, the central limit theorem
\[
Z_n= \sum_{i=1}^{k_n} \wt
S_n(K_{ni}) \std N\bigl(0,\bfc'
\Sigma_N\bfc\bigr),
\]
holds if and only if the following three conditions are satisfied:
$EZ_n=0$, $\var(Z_n)\to\bfc'\Sigma_N\bfc$ and for every $\vep>0$,
%
\begin{equation}
\label{eqq} \sum_{i=1}^{k_n} E \bigl[ \bigl(
S_n (K_{ni})\bigr)^2 I_{\{ |S(K_{ni})|>
\varepsilon\}} \bigr]
\to0.
\end{equation}
The condition $EZ_n=0$ holds since $E\wt I_t=0$, hence
$E \wt S_n(K_{ni}) =0$ for every $i$.
As for (6.8) in Davis and
Mikosch
\cite{davismikosch2009}, a trivial bound of the left-hand side in
\eqref{eqq} is given by
\[
c \frac{m_n^3}{n} \sum_{i=1}^{k_n} P
\bigl(\bigl|S_n(K_{ni})\bigr|> \varepsilon\bigr)\le c
\frac{m_n^3}n \sum_{i=1}^{k_n}I_{\{ c (m_n^3/n)^{0.5}>
\varepsilon\}}
.
\]
In view of (M1), $m_n^3/n=\mathrm{o}(1)$, and therefore the right-hand side vanishes for
sufficiently large~$n$. Therefore, \eqref{eqq} holds.
%
\begin{lemma}Under the conditions of Theorem~\ref{thm1},
\[
\var(Z_n)=\sum_{i=1}^{k_n} \var
\bigl(S_n(K_{ni})\bigr)\to\bfc'
\Sigma_N\bfc.
\]
\end{lemma}
\begin{pf}
We proceed in a similar way as for Proposition \ref{prop1}.
It will be convenient to introduce the following notation for $\la\in
(0,\uppi)$,
\begin{eqnarray*}
\wt\alpha_n(\la)&=& \biggl(\frac{2m_n}n \biggr)^{1/2}
\sum_{i=1}^{k_n} \sum
_{t\in K_{ni}} \cos(\la t) \wt I_t(i) ,
\\
\wt\beta_n(\la)&=& \biggl(\frac{2m_n}n \biggr)^{1/2}
\sum_{i=1}^{k_n} \sum
_{t\in K_{ni}} \sin(\la t) \wt I_t(i) ,
\end{eqnarray*}
where for each $i\le k_n$,
\[
(I_{1},\ldots,I_{m_n})\eqd\bigl(I_{(i-1)m_n+1}(i),
\ldots,I_{im_n}(i)\bigr)
\]
the vectors on the right-hand side are mutually independent for $i\le
k_n$ and the
quantities $\wt I_t(i)$ are the mean corrected versions of
$I_t(i)$, that is, $\wt I_t(i)=I_t(i)-p_0$. The statement of the lemma
is proved if we can show that the pairs
$(\wt\alpha_n(\la),\wt\beta_n(\w)$,
$(\wt\alpha_n(\la),\wt\alpha_n(\w)$,
$(\wt\beta_n(\la),\wt\beta_n(\w)$, $(\wt\alpha_n(\la),\wt
\beta_n(\la))$,
are asymptotically uncorrelated for $\la\ne\w$ and that
%
\begin{equation}
\label{eqqq} \var\bigl(\wt\alpha_n(\la)\bigr)\sim\var\bigl(\wt
\beta_n(\la)\bigr)\sim \mu_0(A) \Biggl[1+2 \sum
_{h=1}^\infty\rho_A(h)\cos(\la h)
\Biggr].
\end{equation}
We check the asymptotic variance of $\wt\alpha_n(\la)$ and omit
similar calculations for $\var(\wt\beta_n(\la))$.
By independence of the sums over the blocks $K_{ni}$ we have for fixed
$k\ge1$,
\begin{eqnarray*}
&&\var\bigl(\wt\alpha_n(\la)\bigr)\\
&&\quad=2\frac{m_n}{n} \sum
_{i=1}^{k_n} \var \biggl(\sum
_{t\in K_{ni}} \cos(\la t) \wt I_t \biggr)
\\
&&\quad=2\frac{m_n}{n} \Biggl[\sum_{i=1}^{k_n}
\sum_{t\in K_{ni}} \var(I_t) \bigl(\cos(\la t)
\bigr)^2 + \sum_{i=1}^{k_n}\sum
_{(i-1)m_n+1\le t\ne s\le im_n}\cov(I_t,I_s) \cos(
\la t)\cos(\la s) \Biggr]
\\
&& \quad =2\frac{m_n}{n} \bigl(p_0-p_0^2
\bigr) \sum_{t=1}^{n} \bigl(\cos(\lambda t)
\bigr)^2
\\
&&\qquad{} +2\frac{m_n}{n} \sum_{i=1}^{k_n}
\sum_{h=1}^{m_n-1} \sum
_{t=1}^{m_n-h} \bigl(p_h-p_0^2
\bigr) (\cos(\lambda h)+\cos \bigl( \la h+2\la\bigl( t+ (i-1) m_n
\bigr) \bigr)
\\
&&\quad=P_1+P_{21}+P_{22}.
\end{eqnarray*}
Then we have by (M) and regular variation of $(X_t)$,
\[
P_1+P_{21} \sim\mu_0(A)+2 \sum
_{h=1}^{m_n-1} \bigl(p_h-p_0^2
\bigr) (m_n-h)\cos(\lambda h)\sim \mu_0(A)f_A(
\la).
\]
We have for fixed $k\ge1$,
\[
2\frac{m_n}{n} \Biggl| \sum_{i=1}^{k_n}\sum
_{h=k+1}^{m_n-1} \sum
_{t=1}^{m_n -h} \bigl(p_{h}-p_0^2
\bigr) \cos\bigl( \la h +2\la\bigl(t+ (i-1) m_n\bigr) \bigr) \Biggr|
\leq c m_n \sum_{h=k+1}^{m_n -1}
\bigl|p_{h}-p_0^2\bigr|,
\]
and the right-hand side is negligible in view of (M) by first letting
$n\to\infty$ and then $k\to\infty$.
Thus, it suffices to consider only finitely many $h$-terms in
$P_{22}$. In view of \eqref{eqa}, for fixed $k$ as $\nto$,
\[
\Biggl\llvert 2\frac{m_n}{n} \sum_{i=1}^{k_n}
\sum_{h=1}^k \bigl(p_h-p_0^2
\bigr) \sum_{t=1}^{m_n -h} \cos\bigl( \la h+2\la
\bigl( t+(i-1) m_n\bigr) \bigr) \Biggr\rrvert \leq c\sum
_{h=1}^k \bigl|p_{h}-p_0^2\bigr|=\mathrm{o}(1)
.
\]
This proves \eqref{eqqq}.

Next, we consider the case of two different frequencies $\la,\w\in
(0,\uppi)$ and show that the following covariances vanish as $\nto$:
\begin{eqnarray*}
&&\cov\bigl(\wt\alpha_n(\la),\wt\alpha_n(\w)
\bigr)
\\
&&\quad= \frac{2m_n}{n} \sum_{i=1}^{k_n}
\cov \Biggl( \sum_{t=1}^{m_n} \wt
I_t \cos\bigl(\lambda \bigl(t+(i-1)m_n\bigr)\bigr),
\sum_{t=1}^{m_n} \wt I_t \cos
\bigl(\omega\bigl(t+(i-1)m_n\bigr)\bigr) \Biggr)
\\
& &\quad=\frac{2m_n}{n}\sum_{t=1}^{n}
\bigl(p_0 -p_0^2\bigr) \cos(\lambda t)\cos
(\omega t)
\\
&&\qquad{} + \frac{2m_n}{n} \sum_{i=1}^{k_n}
\sum_{h=1}^{m_n -1} \sum
_{t=1}^{m_n -h} \bigl(p_h -
p_0^2\bigr) \bigl[\cos\bigl(\lambda
\bigl(t+(i-1)m_n+h\bigr)\bigr) \cos(\omega\bigl(t+(i-1)m_n
\bigr)
\\
& &\hspace*{133pt}\qquad{} +\cos\bigl(\lambda\bigl(t+(i-1)m_n\bigr)\bigr)\cos\bigl(\omega
\bigl(t+(i-1)m_n +h\bigr)\bigr) \bigr]\\
&&\quad =Q_1+Q_2
.
\end{eqnarray*}
In view of \eqref{eqa} and since $\la\ne\w$,
\begin{eqnarray*}
|Q_1|&= & \frac{m_n}{n} \bigl(p_0-p_0^2
\bigr) \Biggl| \sum_{t=1}^{n} \bigl(\cos \bigl((
\lambda +\omega)t\bigr) +\cos\bigl((\lambda-\omega)t\bigr) \bigr) \Biggr| \le c
\frac{m_n}{n} \bigl(p_0 -p_0^2
\bigr)=\mathrm{O}\bigl(n^{-1}\bigr). 
\end{eqnarray*}
Similarly, multiple application of \eqref{eqa}, first summing over
$t$, then over $l$, yields
\begin{eqnarray*}
|Q_2|&=& \frac{m_n}{n} \Biggl| \sum_{h=1}^{m_n}
\bigl(p_h -p_0^2\bigr) \sum
_{l=0}^{k_n -1} \sum_{t=1}^{m_n -h}
\bigl( \cos\bigl((\lambda+ \omega) (t +h +lm_n)+\lambda h \bigr)
\\
& & \hspace*{116pt}{}+\cos\bigl((\lambda-\omega) (t +h + l m_n)+ \lambda h\bigr)\\
&&\hspace*{116pt}{} +
\cos \bigl((\lambda+ \omega) (t +h +lm_n)+\omega h \bigr)
\\
& &\hspace*{116pt} {}+ \cos\bigl((\lambda- \omega) (t +h +lm_n)- \omega h \bigr) \bigr)\Biggr |
\\
& \leq&
c_0 \sum_{h=1}^{m_n}
\bigl|p_h-p_0^2\bigr| \leq c \frac{r_n}{m_n}
(m_np_0) + c \frac{r_n}{m_n^2} (m_np_0)^2
+ c\sum_{h=r_n+1}^{m_n}\xi _h \to0 ,
\end{eqnarray*}
where $c_0= 4\max\{1/\sin((\lambda+\omega)/2),1/\sin(|\lambda
-\omega|/2)\}+4$.
Thus $\cov(\wt
\alpha_n(\la),\wt\alpha_n(\w))=\mathrm{o}(1)$. Using similar arguments, it also
follows that the covariances of the pairs $(\wt\alpha_n(\la),\wt
\beta_n(\w))$,
$(\wt\beta_n(\la),\wt\beta_n(\w))$ and $(\wt\alpha_n(\la),\wt
\beta_n(\la))$ are asymptotically negligible. This proves the lemma.
\end{pf}

\section{Smoothing the periodogram}\label{subsecsmooth}
Corollary~\ref{cor1} is analogous to the asymptotic theory for
the periodogram of a stationary sequence; see Brockwell and
Davis \cite{brockwelldavis1991}, Section~10.4, where the
corresponding results are
proved for the periodogram ordinates of a general linear processes
with i.i.d. innovations. These results are then employed for
showing that smoothed versions of the periodogram are consistent
estimators of the spectral density at a given frequency.
Our next goal is to prove a similar result.

We start by introducing the smoothed periodogram. For a fixed
frequency $\la\in(0,\uppi)$ define
\begin{eqnarray*}
\la_{0} = \min\{2\uppi j/n \dvt 2\uppi j/n \geq \lambda\} \quad \mbox{and}\quad
\lambda_{j} = \lambda_{0} +2\uppi j/n ,\qquad |j|\le
s_n.
\end{eqnarray*}
Here we suppress the dependence of $\la_j$ on $n$.
In what follows, we will assume that $s_n\to\infty$ and $s_n/n\to0$
as $\nto$. For a given set $A\subset\ov\bbr_0^d$ bounded away from
zero and
any non-negative weight function  $w=(w_n(j)) _{|j|\le s_n}$
satisfying the
conditions
%
\begin{equation}
\label{eqweight}  \sum_{|j|\leq s_n} w_n (j) =1
\quad\mbox{and}\quad \sum_{|j|\leq s_n} w_n^2
(j) \to0 \qquad\mbox{as } n \to\infty,
\end{equation}
we introduce the \emph{smoothed periodogram}
\[
\widetilde{f}_{nA} (\lambda) = \sum_{|j| \le s_n}
w_n(j) I_{nA} (\lambda_{j}).
\]
%
\begin{theorem}\label{thmsmoth}
Assume the conditions of Theorem~\ref{thm1}, \eqref{eqweight} on the
weight function $w$ and \textup{(M2)}.
Then for every fixed frequency $\la\in(0,\uppi)$, as $\nto$,
\begin{eqnarray*}
\wt f_{nA}(\la)\stackrel{L^2} {\rightarrow}
\mu_0(A) f_A(\la)\quad \mbox{and}\quad \frac{ \wt f_{nA}(\la)}{\wh P_m(A)}\stp
f_A(\la).
\end{eqnarray*}
\end{theorem}

\begin{figure}

\includegraphics{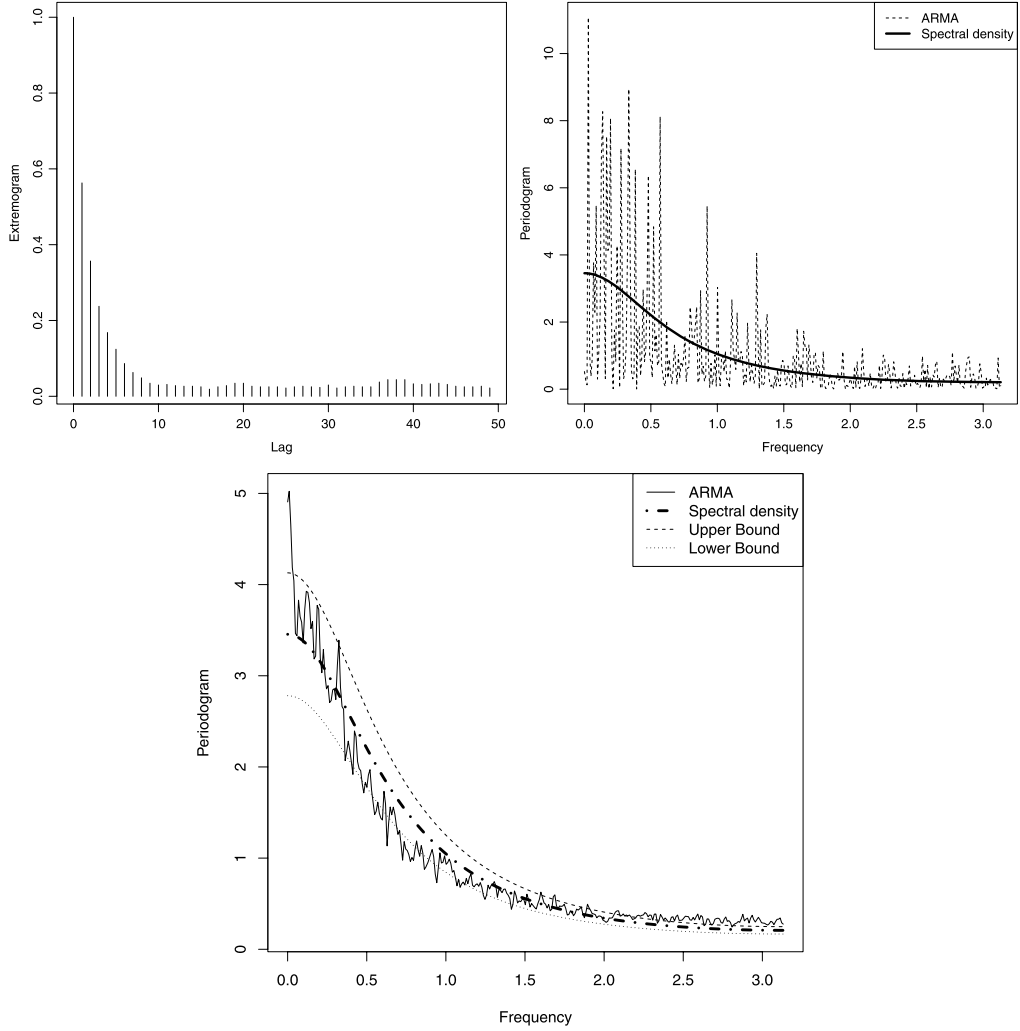}

\caption{Top-Left: The sample extremogram of an $\operatorname{ARMA}(1,1)$
process with parameters $\phi=0.8$, $\theta=0.1$ and
i.i.d. $t$-distributed noise with 3 degrees of freedom. We choose
$A=(1,\infty)$.
Top-Right: The corresponding raw periodogram and the theoretical
spectral density
$f_A$ (solid line). Bottom: The smoothed
periodogram with Daniell window, $s_n=50$.}\label{fig1}
\end{figure}

\begin{figure}

\includegraphics{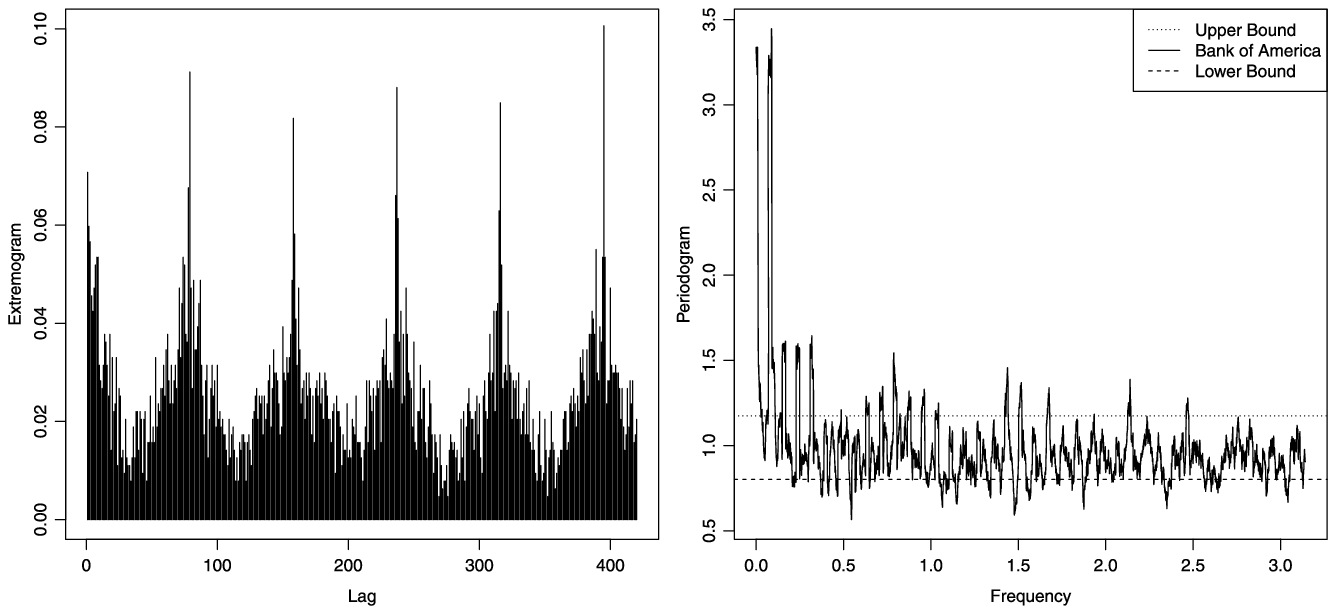}

\caption{Left: The sample extremogram
of 5-min returns of BAC stock price for $A=(1,\infty)$. Right:
The smoothed periodogram with Daniell window, $s_n=50$. The confidence
bands are constructed from the smoothed periodograms of 99 permutations of
the data.}
\label{fig2}
\end{figure}

In Figures~\ref{fig1} and \ref{fig2} we show the extremogram, the
standardized periodogram and the corresponding smoothed periodogram for some simulated and
real-life data.
The data underlying Figure~\ref{fig1} are simulated from an $\operatorname{ARMA}(1,1)$ process
$(X_t)$ with parameters $\phi=0.8$ and $\theta=0.1$
and i.i.d. $t$-distributed noise $(Z_t)$ with 3 degrees of freedom,
hence $(X_t)$ is
regularly varying with $\alpha=3$.
The top-left graph shows the sample extremogram based
on a sample of size $n=31\mbox{,}757$ and the threshold is chosen as the $98\%
$ empirical quantile of the data.
The top-right graph visualizes the theoretical spectral density $f_A$
for $A=(1,\infty)$ (see Appendix~\ref{exam1} for an expression)
and the raw periodogram which exhibits rather erratic behavior. The bottom
graph shows the smoothed periodogram with Daniell window
$w_n(i)=1/(2 s_n+1)$, $|i|\le s_n=50$. We also show the curves
$f_A(\lambda) (1 \pm1.96/\sqrt{2 s_n
+1})$, which constitute a confidence band based on the following heuristic
argument. In the proof of Theorem \ref{thmsmoth}, we show that $\var
(\widetilde{f}_{nA} (\lambda) ) \sim\sum_{|j|\le s_n } w_n^2(j)\mu
_0^2(A) f_{nA}^2
(\lambda)$ for every $\lambda\in(0,\uppi)$. Furthermore, we know that
$\widehat{P}_m (A)
\stackrel{P}{\to} \mu_0(A)$. Based on these calculations, we take
$\sum_{|j|\le s_n} w_n^2(j) f_{nA}^2(\lambda)$ as a surrogate
quantity for the unknown variance of
$\widetilde{f}_n (\lambda)/\widehat{P}_m (A)$.

The data underlying Figure~\ref{fig2} are 5-min returns for
the stock price of Bank of America (BAC) with the sample size
$n=31\mbox{,}757$, and
$a_m$ is chosen as the $98\%$ empirical quantile of the data. We
provide the same type of analysis as in Figure~\ref{fig1} for these data.
The largest peak in the periodogram at the frequency $0.29$ corresponds
to an extremal cycle length of
$6$ hours, this is roughly the length of a trading day. We also show
$95\%$ pointwise confidence bands for
the smoothed periodogram. They are not asymptotic since we do not have
a central limit theorem for the smoothed periodogram yet.
They are constructed from the distribution of the corresponding
smoothed periodogram
s based on 99 random permutations of the
data. If the data were i.i.d., any permutation would not change
the dependence
structure of the data and one would expect that the estimated spectral
density stays inside the band, but this is obviously not the case,
indicating that the data exhibit some significant extremal dependence.

%

\begin{pf*}{Proof of Theorem \ref{thmsmoth}}
We mentioned in Remark~\ref{rem2} that
%
\begin{equation}
\label{eqlll} EI_{nA} (\lambda)\to\mu_0(A)
f_A(\la) \qquad\mbox{as $\nto$ uniformly on sets $[a,b]\subset (0,\uppi)$.}
\end{equation}
Therefore, since $\max_{|j|\le s_n}|\la_{j}-\la|\to0$
and $f_A$ is continuous, we have
\[
E\widetilde{f}_{nA} (\lambda)=\sum_{|j|\le s_n}w_n(j)
EI_{nA} (\lambda_{j}) \to\mu_0(A)
f_A(\la) , \qquad\nto.
\]
The statement of the theorem then follows if we can show that
$\var( \widetilde{f}_n (\lambda))\to0$. We observe that
\begin{eqnarray*}
\var \bigl( \widetilde{f}_{nA} (\lambda) \bigr) = \sum
_{|j|\leq s_n} w_n^2(j) c_{jj} +
\sum_{-s_n\le j_1\ne j_2 \leq s_n
} w_n(j_1)
w_n(j_2) c_{j_1j_2}.
\end{eqnarray*}
In view of condition \eqref{eqweight} it suffices to show that
$c_{j_1j_2}=\cov(I_{nA} (\lambda_{j_1}) , I_{nA} (\lambda_{j_2}))\to0$
and
%
\begin{eqnarray}
\label{eq14a} c_{jj}=\var\bigl(I_{nA} (
\lambda_{j})\bigr)\to\bigl(\mu_0(A) f_A(
\lambda )\bigr)^2\qquad \mbox{uniformly for $j,j_1,j_2
\in[-s_n, s_n]$, $j_1\ne j_2$.}\qquad
\end{eqnarray}
We will only show \eqref{eq14a}; the proof of $c_{j_1,j_2}\to0$ for
$j_1\ne j_2$ is similar and therefore omitted.
Since \eqref{eqlll} holds, we have to show that
%
\begin{equation}
\label{eqsuff} E \bigl(I_{nA}^2(\lambda_{j})
\bigr) \to2 \bigl(\mu_0(A) f_A(\lambda)
\bigr)^2.
\end{equation}
Recall $\wh f_{nA}(\la)$ from \eqref{eqlagwindow}
and define
\[
\wh g_{nA}(\la)=2 \sum_{h=r_n+1}^{n-1}
\cos(\lambda h)\widetilde {\gamma}_n(h).
\]
We will study the decomposition
\[
E\bigl(I_{nA}^2(\lambda_{j})\bigr) = E \wh
f_{nA}^2(\la_{j}) + 2 E \bigl(\wh
f_{nA}(\la_{j}) \wh g_{nA}(\la_{j})
\bigr) + E\wh g_{nA}^2(\la_{j}).
\]
Following the lines of the proof of Theorem 5.1
in \cite{davismikosch2009}, we
conclude that
%
\begin{equation}
\label{eq101} E \wh f_{nA}^2(\la_{j}) \to
\bigl(\mu_0(A) f_A(\lambda)\bigr)^2 ,
\end{equation}
uniformly for the considered frequencies $\la_{j}$. Then
\eqref{eqsuff} is proved if we can show that
%
\begin{eqnarray}
\label{eqhe1} E\bigl(\wh f_{nA}(\la_{j}) \wh
g_{nA}(\la_{j})\bigr)&\to& 0 ,
\\
E\wh g_{nA}^2(\la_{j})&\to& \bigl(
\mu_0(A) f_A(\lambda)\bigr)^2.
\label{eqhe2}
\end{eqnarray}
Throughout we will use the notation, for $h_1,h_2,h_3\ge0$,
\begin{eqnarray*}
p_{h_1h_2h_3}&=&P(X_0>a_m,X_{h_1}>a_m,X_{h_1+h_2}>a_m,X_{h_1+h_2+h_3}>a_m)
,
\\
p_{h_1h_2}&=&p_{h_1h_20} ,\qquad p_{h_1}=p_{h_10} ,
\end{eqnarray*}
and we observe that
%
\begin{eqnarray}
\label{eqhelp1} p_h&=&\bigl(p_h-p_0^2
\bigr)+p_0^2 ,
\\
p_{h_1h_2}&=& (p_{h_1h_2}-p_{h_1}p_0)+p_{h_1}p_0=p_{h_1h_2}-p_0p_{h_2}+p_0p_{h_2}
\label {eqhelp2}
\nonumber
\\[-8pt]
\\[-8pt]
\nonumber
&=&(p_{h_1h_2}-p_0p_{h_2})+p_0
\bigl(p_{h_2}-p_0^2\bigr)+p_0^3
,
\\
p_{h_1h_2h_3}&=&(p_{h_1h_2h_3}-p_0p_{h_2h_3})+p_0p_{h_2h_3}
\label {eqhelp3}
\nonumber
\\[-8pt]
\\[-8pt]
\nonumber
&=& (p_{h_1h_2h_3}-p_0p_{h_2h_3})+p_0(p_{h_2h_3}-p_0p_{h_3})+p_0^2p_{h_3}
.
\nonumber
\end{eqnarray}
\subsection*{Proof of \protect\eqref{eqhe1}}
We have
\begin{eqnarray*}
E\bigl(\wh f_{nA}(\la_{j}) \wh g_{nA}(
\la_{j})\bigr) &=&E \Biggl[2\widetilde{\gamma}_n(0) \wh
g_{nA}(\la_j) +4 \wh g_{nA}(
\la_j) \sum_{h=1}^{r_n}\cos(
\la_j h)\wt\gamma_n(h) \Biggr]
\\
&=&
J_{1}+J_{2} ,
\end{eqnarray*}
where
\begin{eqnarray*}
J_{1}&=& 4 \frac{m_n^2}{n^2} \sum_{t_1=1}^n
\sum_{h=r_n+1}^{n-1}\sum
_{t_2=1}^{n-h}E[ I_{t_1}I_{t_2}I_{t_2+h}]
\cos(\lambda_{j} h) ,
\\
J_2&=&8 \frac{m_n^2}{n^2}\sum_{t_1=1}^{n-1}
\sum_{h_1=1}^{r_n}\sum
_{h_2=r_n+1}^{n-1}\sum_{t_2=1}^{n-h_2}E[I_{t_1}I_{t_1+h_1}I_{t_2}I_{t_2+h_2}]
\cos(\lambda _{j}h_1)\cos(\lambda_{j}h_2)
.
\end{eqnarray*}

\subsubsection*{Proof that $J_1$ is negligible}
We observe, that depending on the values $h, t_1,t_2$,
$E[I_{t_1}I_{t_2}I_{t_2+h}]$ may simplify: if $t_1=t_2$ or
$t_1=t_2+h$, $E[ I_{t_1}I_{t_2}I_{t_2+h}]=p_h$; if $t_1<t_2$, $E[
I_{t_1}I_{t_2}I_{t_2+h}]=p_{t_2-t_1,h}$;
if $t_2<t_1<t_2+h$, $E[
I_{t_1}I_{t_2}I_{t_2+h}]=p_{t_1-t_2,h-t_1+t_2}$; if $t_1>t_2+h$, $E[
I_{t_1}I_{t_2}I_{t_2+h}]=p_{h,t_1-h-t_2}$.
If we take into account these different cases, we obtain
\begin{eqnarray*}
J_1&=&4 \frac{m_n^2}{n^2}\sum_{h=r_n+1}^{n-1}(n-h)
(2p_h)\cos(\lambda _{j}h)+4 \frac{m_n^2}{n^2}\sum
_{h_2=r_n+1}^{n-2}\sum
_{h_1=1}^{n-h_2-1}(n-h_1-h_2)p_{h_1h_2}
\cos(\lambda_{j}h_2)
\\
&&{} +4 \frac{m_n^2}{n^2}\sum_{h_2=r_n+1}^{n-1}
\sum_{h_1=1}^{h_2-1}(n-h_2)p_{h_1,h_2-h_1}
\cos(\lambda_{j}h_2)
\\
&&{} +4 \frac{m_n^2}{n^2}\sum_{h_2=r_n+1}^{n-2}
\sum_{h_1=1}^{n-h_2-1}(n-h_1-h_2)p_{h_2h_1}
\cos (\lambda_{j}h_2)\\
&=&\sum_{i=1}^4
J_{1i}.
\end{eqnarray*}
Applying \eqref{eqhelp1}, the mixing condition (M2) and
Lemma~\ref{lemaa} imply that
\[
J_{11}\le c m_n\sum_{h=r_n+1}^\infty
\xi_h + c \frac{(m_np_0)^2}{n(\sin(\lambda_j /2))^2}= \mathrm{o}(1).
\]
As regards $J_{12}$, apply \eqref{eqhelp2} and split the $h_1$-index
set into $h_1\le r_n$ and $h_1>r_n$. Then (M2) and Lemma~\ref{lemaa}
imply that
\begin{eqnarray*}
|J_{12}| &\le& 
c m_n^2\sum
_{h_2=r_n+1}^{n-1}\xi_{h_2}
\\
&&{}+c \Biggl| \frac{m_n^2}{n^2}\sum_{h_2=r_n+1}^{n-2}
\Biggl(\sum_{h_1=1}^{\min
(r_n,n-h_2-1)} +\sum
_{h_1=r_n+1}^{n-h_2-1} \Biggr) (n-h_1-h_2)
\bigl(p_{h_1}\pm p_0^2\bigr)p_0
\cos(\lambda_{j}h_2)\Biggr |
\\
&\le& \mathrm{o}(1)+ c\frac{r_n}{n}(m_np_0)^2
\bigl(\sin(\lambda_{j}/2)\bigr)^{-2}+c (m_np_0)m_n
\sum_{h_1=r_n+1}^{n-1}\xi_{h_1}+c
\frac
{(m_np_0)^3}{m_n}=\mathrm{o}(1).
\end{eqnarray*}
Now consider $J_{13}$. Abusing notation, we will write $h_2$ instead of
$h_2-h_1$. Introduce the index sets
\begin{eqnarray*}
K_1&=&\bigl\{(h_1,h_2)\dvt 1\le
h_i\le r_n ,i=1,2\bigr\} ,
\\
K_2&=&\bigl\{(h_1,h_2)\dvt 1\le
h_1\le r_n ,r_n<h_2<n-h_1
\bigr\} ,
\\
K_3&=&\bigl\{(h_1,h_2)\dvt
r_n< h_1\le n-1 ,1\le h_2\le
\min(r_n,n-h_1-1)\bigr\} ,
\\
K_4&=&\bigl\{(h_1,h_2)\dvt
r_n< h_1\le n-1 ,r_n< h_2<
n-h_1\bigr\}.
\end{eqnarray*}
Now introduce the mixing coefficients $\xi_h$ and use Lemma~\ref{lemaa}:
\begin{eqnarray*}
|J_{13}| &\le& c \frac{m_n^2}{n^2}\Biggl |\sum
_{h_1=1}^{n-2}\sum_{h_2=\max
(1,r_n+1-h_1)}^{n-h_1-1}(n-h_1-h_2)p_{h_1h_2}
\cos\bigl(\lambda _{j}(h_1+h_2)\bigr) \Biggr|
\\
&\leq& c \frac{m_n^2}{n^2}\sum_{i=1}^4
\Biggl|\sum_{K_i} (n-h_1-h_2)p_{h_1h_2}
\cos\bigl(\lambda_{j}(h_1+h_2)\bigr)\Biggr |
\\
&\leq& c \frac{m_nr_n^2}{n} (m_np_0) + c \Biggl[
\frac{m_nr_n}{n} m_n \sum_{h_2=r_n+1}^{n-1}
\xi_{h_2}+\frac
{r_n}{n}(m_np_0)^2
\bigl(\sin(\lambda_{j}/2)\bigr)^{-2} \Biggr]
\\
& &{} +c \Biggl[ \frac{m_nr_n}{n} m_n \sum
_{h_1=r_n+1}^{n-1}\xi_{h_1}+ \frac{r_n}{n}(m_np_0)^2
\bigl(\sin(\lambda_{j}/2)\bigr)^{-2} \Biggr]
\\
& &{} +c \Biggl[ m_n^2 \sum_{h_1=r_n+1}^{n-1}
\xi_{h_1}+(m_np_0)m_n\sum
_{h_2=r_n+1}^{n-1}\xi_{h_2}+\frac{1}{m_n}(m_np_0)^3
\bigl(\sin(\lambda _{j}/2)\bigr)^{-2} \Biggr].
\end{eqnarray*}
The right-hand side vanishes as $\nto$ by virtue of (M2). The same
idea of proof
applies to the relation $J_{14}=\mathrm{o}(1)$. Thus, we showed that $J_1=\mathrm{o}(1)$.

\subsubsection*{Proof that $J_2$ is negligible} We split the
summation over
disjoint index sets, depending on the ordering of
$\{t_1,t_1+h_1,t_2,t_2+h_2 \}$:
$t_1=t_2$, $t_1+h_1=t_2+h_2$, $t_1+h_1=t_2$, $t_1=t_2+h_2$,
$t_1<t_2<t_1+h_1<t_2+h_2$, $t_2<t_1<t_1+h_1<t_2+h_2$,
$t_2<t_1<t_2+h_2<t_1+h_1$, $t_1+h_1<t_2$ and $t_2+h_2<t_1$.
Consider the index sets (we recycle the notation $h_1,h_2$ here)
\begin{eqnarray*}
L_1&=&\bigl\{(h_1,h_2)\dvt 1\le
h_1\le r_n , r_n<h_2<n\bigr\} ,
\\
L_2&=&\bigl\{(h_1,h_2)\dvt 1\le
h_1\le r_n , r_n<h_2<n-h_1
\bigr\} ,
\\
L_3&=&\bigl\{(h_1,h_2,h_3)\dvt
2\le h_1\le r_n , r_n<h_2<n-h_1-1
,1\le h_3<h_1\bigr\} ,
\\
L_4&=&\bigl\{(h_1,h_2,h_3)
\dvt 1\le h_1\le r_n , r_n<h_2<n,
1\le h_3<h_2\bigr\} ,
\\
L_5&=&\bigl\{(h_1,h_2,h_3)
\dvt 1\le h_1\le r_n , r_n<h_2<n-1,
h_2-h_1< h_3\le\min(n,h_2+h_1-1)
\bigr\} ,
\\
L_6&=&\bigl\{(h_1,h_2,h_3)
\dvt 1\le h_1\le r_n , r_n<h_2<n-h_1-1,
1\le h_3< n-h_1-h_2\bigr\}.
\end{eqnarray*}
We write for short $f_{h_1h_2}=\cos(\lambda_{j_1}h_1) \cos(\lambda
_{j_1}h_2)$.
Then
\begin{eqnarray*}
J_{2}&= & 8\frac{m_n^2}{n^2} \biggl[ \sum
_{L_1}(n-h_2) (p_{h_1,h_2-h_1}+p_{h_2-h_1,h_1}
) f_{h_1h_2}+ \sum_{L_2} (n-h_1-h_2)
(p_{h_2h_1}+p_{h_1h_2} ) f_{h_1h_2}
\\
& &\hspace*{23pt}{}+ \sum_{L_3} (n-h_2-h_3)
p_{h_3,h_1-h_3,h_2-h_1+h_3} f_{h_1h_2} +\sum_{L_4}
(n-h_2) p_{h_3,h_1,h_2-h_1} f_{h_1h_2}
\\
& &\hspace*{23pt}{}+\sum_{L_5} (n-h_1-h_3)
p_{h_3,h_2-h_3,h_1-h_2+h_3} f_{h_1h_2}
\\
&& \hspace*{23pt}{}+\sum_{L_6}(n-h_1-h_2-h_3)
(p_{h_1h_3h_2}+p_{h_2h_3h_1} ) f_{h_1h_2} \biggr] 
\\
&=&\sum_{i=1}^6
J_{2i}.
\end{eqnarray*}
The terms $J_{2i}$, $i=1,2$, involve probabilities of the form $p_{kl}$.
These terms can
be treated in the same way as $J_1$ and shown to be negligible. We
omit details.

The remaining $J_{2i}$'s contain probabilities of the form
$p_{kls}$. We illustrate how one can deal with these pieces. We start
with
\begin{eqnarray*}
|J_{23}|& =& 8 \Biggl| \frac{m_n^2}{n^2} \sum
_{h_1=1}^{r_n-1}\sum_{h_3=1}^{r_n-h_1}
\Biggl( \sum_{h_2=r_n+1-h_3}^{r_n}+ \sum
_{h_2=r_n+1}^{n-h_1-h_3-1} \Biggr) (n-h_1-h_2-h_3)p_{h_1h_3h_2}f_{h_1+h_3,h_2+h_3}
\Biggr|
\\
&\le& c \frac{m_n r_n^3}{n} (m_np_0)+ \Biggl[c
\frac{m_n^2}{n} \sum_{h_1=1}^{r_n-1}\sum
_{h_3=1}^{r_n-h_1} \sum
_{h_2=r_n+1}^{n-h_1-h_3-1}|p_{h_1h_3h_2}-p_{h_1h_3}p_0|
\\
&&\hspace*{76pt}{}+c \frac{m_n^2}{n^2}p_0 \sum_{h_1=1}^{r_n-1}
\sum_{h_3=1}^{r_n-h_1}p_{h_1h_3}\cos\bigl(
\lambda_{j_1}(h_1+h_3)\bigr)
\\
&&\hspace*{76pt}\phantom{+}{}\times\sum_{h_2=r_n+1}^{n-h_1-h_3-1}
(n-h_1-h_2-h_3)\cos\bigl(\lambda
_{j_1}(h_1+h_2)\bigr) \Biggr]
\\
&\le& c \frac{m_n r_n^3}{n}+ c \frac{m_n^2 r_n^2}{n} \sum
_{h_2=r_n+1}^\infty\xi_{h_2}+c \frac{r_n^2}{n}
(m_n p_0)^2 \bigl(\sin(
\lambda_j/2)\bigr)^{-2}.
\end{eqnarray*}
In the last step we used Lemma~\ref{lemaa}. The right-hand side in
the latter relation
converges to zero in view of the assumptions on $r_n,m_n$ and (M2).
The remaining expressions $J_{2i}$ which contain probabilities $p_{kls}$
over index sets such that $k,l>r_n,s\le r_n$ or $k>r_n, l,s\le r_n$ can
be shown to be negligible by using similar arguments. We omit details.
Those sums which contain probabilities $p_{kls}$
over index sets such that $k,l,s>r_n$ are most difficult to deal with.
The corresponding bounds follow from the next lemma.
%
\begin{lemma}Let $\la,\w\in[a,b]$, $0<a<b<\uppi$, possibly depending
on $n$,
and $x_1$, $x_2$ be real numbers.
Assume that
%
\begin{equation}
\label{eqcc} m_n^2n \sum_{h=r_n+1}^n
\xi_h \to0 , \qquad \nto,
\end{equation}
where $(\xi_t)$ is the mixing rate function.
Then
%
\begin{eqnarray}
Q_0 = \frac{m_n^2}{n^2}\sum_{h_1,h_2,h_3>r_n}
(n-h_1-h_2-h_3)_+p_{h_1h_2h_3}\cos(
\lambda h_1+x_1)\cos(\omega h_3+x_2)
&\to&0 ,\qquad\label{eqhelp4}
\\
 \frac{m_n^2}{n^2}\sum_{h_1,h_2,h_3>r_n}
(n-h_1-h_2-h_3)_+p_{h_1h_2h_3}\sin(
\lambda h_1+x_1)\sin(\lambda h_3+x_2)
&\to&0. \label{eqhelp5}
\end{eqnarray}
\end{lemma}
\begin{pf}
Recall \eqref{eqhelp3}. Write $g_{h_1h_3}=\cos(\lambda h_1+x_1)\cos
(\omega h_3+x_2)$.
Then we have
\begin{eqnarray*}
|Q_0| &\leq& \frac{m_n^2}{n^2}\sum_{h_1,h_2,h_3>r_n}
(n-h_1-h_2-h_3)_+|p_{h_1h_2h_3}-p_{h_1}p_{h_3}|
\\
&&{}+ \Biggl|\frac{m_n^2}{n^2}\sum_{h_3=r_n+1}^{n-2r_n-3}\sum
_{h_1=r_n+1}^{n-h_3-r_n-2}\sum
_{h_2=r_n+1}^{n-h_1-h_3-1}(n-h_1-h_2-h_3)
\bigl(p_{h_1}-p_0^2\bigr)
\bigl(p_{h_3}-p_0^2\bigr) g_{h_1h_3} \Biggr|
\\
&&{}+ \Biggl|\frac{m_n^2}{n^2}\sum_{h_2=r_n+1}^{n-2r_n-3}\sum
_{h_3=r_n+1}^{n-h_2-r_n-2}\sum
_{h_1=r_n+1}^{n-h_2-h_3-1}(n-h_1-h_2-h_3)p_0^2
\bigl(p_{h_3}-p_0^2\bigr)g_{h_1h_3} \Biggr|
\\
&&{}+ \Biggl|\frac{m_n^2}{n^2}\sum_{h_2=r_n+1}^{n-2r_n-3}\sum
_{h_1=r_n+1}^{n-h_2-r_n-2}\sum
_{h_3=r_n+1}^{n-h_1-h_2-1}(n-h_1-h_2-h_3)p_0^2
\bigl(p_{h_1}-p_0^2\bigr) g_{h_1h_3} \Biggr|
\\
&&{}+\Biggl |\frac{m_n^2}{n^2}\sum_{h_2=r_n+1}^{n-2r_n-3}\sum
_{h_1=r_n+1}^{n-h_2-r_n-2}\sum
_{h_3=r_n+1}^{n-h_1-h_2-1}(n-h_1-h_2-h_3)p_0^4
g_{h_1h_3} \Biggr|\\
&=& \sum_{i=1}^5
Q_i.
\end{eqnarray*}
By virtue of \eqref{eqcc}, $Q_1$ is negligible. Similarly,
$Q_2\leq m_n^2 (\sum_{h=r_n+1}^n\xi_{h})^2 \to0$.
As to $Q_3$, Lemma~\ref{lemaa} and mixing imply that
\begin{eqnarray*}
Q_3&\leq&c \frac{(m_np_0)^2}{n^2}\sum_{h_2=r_n+1}^{n-2r_n-3}
\sum_{h_3=r_n+1}^{n-h_2-r_n-2}(n-h_2-h_3)\bigl|p_{h_3}-p_0^2\bigr|
\bigl(\sin(\la /2)\bigr)^{-2}\leq c n \sum
_{h_3=r_n+1}^n\xi_{h_3}\to0.
\end{eqnarray*}
A similar bound applies to $Q_4$. A double application of Lemma \ref
{lemaa} yields
\[
Q_5 \leq 
c \frac{(m_np_0)^4}{m_n^2} \bigl(\sin(\omega/2)
\sin(\lambda /2)\bigr)^{-2}\to0.
\]
Collecting these bounds, we proved \eqref{eqhelp4}. Similar arguments
apply to \eqref{eqhelp5}.
\end{pf}
Thus, we showed that $J_1$ and $J_2$ are negligible as $\nto$. Hence,
\eqref{eqhe1} holds.

\subsection*{Proof of \protect\eqref{eqhe2}}
Following the steps for showing that $J_2$ is negligible, we decompose
$E\wh g_{nA}^2(\la_{j})$ into sums over disjoint index sets
depending on the ordering of $\{t_1,t_1+h_1,t_2,t_2+h_2\} $: $t_1=t_2$
and $h_1=h_2$; $t_1=t_2$ and $h_1>h_2$; $t_1=t_2$ and $h_1<h_2$;
$t_1+h_1=t_2+h_2$ and $t_1>t_2$; $t_1+h_1=t_2+h_2$ and $t_1<t_2$;
$t_1=t_2+h_2$; $t_2=t_1+h_1$; $t_1<t_2<t_1+h_1<t_2+h_2$;
$t_2<t_1<t_2+h_2<t_1+h_1$; $t_1<t_2<t_2+h_2<t_1+h_1$;
$t_2<t_1<t_1+h_1<t_2+h_2$; $t_1>t_2+h_2$; $t_2>t_1+h_1$.
Consider the index sets (we recycle the notation $h_1,h_2$ here)
\begin{eqnarray*}
B_1 & =& \{ h\dvt r_n<h<n \} ,
\\
B_2 & =& \bigl\{ (h_1,h_2)\dvt
r_n<h_1<n-r_n, 1\leq h_2<n-h_1
\bigr\} ,
\\
B_3 & =& \bigl\{ (h_1,h_2)\dvt
r_n<h_1<n-r_n-1,r_n<h_2<n-h_1
\bigr\} ,
\\
B_4 & =& \bigl\{ (h_1,h_2,h_3)
\dvt 1\leq h_1 <n-r_n-2, r_n<h_2<n-h_1-1,
1\le h_3<n-h_1-h_2 \bigr\} ,
\\
B_5 & =& \bigl\{ (h_1,h_2,h_3)
\dvt 1\leq h_1 <n-r_n-1, \max(1,r_n+1-h_1)
\leq h_2<n-h_1-1,
\\
&&\hspace*{2pt}{}\max(1,r_n+1-h_2) \leq h_3<n-h_1-h_2
\bigr\} ,
\\
B_6 & =& \bigl\{(h_1,h_2,h_3)
\dvt r_n< h_1 <n-r_n-2, r_n<h_3<n-1-h_1,
1\leq h_2<n-h_1-h_3 \bigr\}.
\end{eqnarray*}
Then we have
\begin{eqnarray*}
E\wh g_{nA}^2(\la_{j}) & = & 4
\frac{m_n^2}{n^2}\sum_{B_1}(n-h)p_hf_{hh}
+4 \frac{m_n^2}{n^2} \sum_{B_2}
(n-h_1-h_2) (p_{h_1h_2}+p_{h_2h_1})f_{h_1+h_2,h_1}
\\
& &{} +4 \frac{m_n^2}{n^2} \sum_{B_3}(n-h_1-h_2)
(p_{h_1h_2}+p_{h_2h_1})f_{h_1h_2}
\\
& &{} +8 \frac{m_n^2}{n^2} (n-h_1-h_2-h_3)
\sum_{B_4} p_{h_1h_2h_3}f_{h_1+h_2+h_3,h_2}
\\
& & {}+8 \frac{m_n^2}{n^2} \sum_{B_5}
(n-h_1-h_2-h_3)p_{h_1h_2h_3}f_{h_1+h_2,h_2+h_3}
\\
&&{} +8 \frac
{m_n^2}{n^2} \sum_{B_6}
(n-h_1-h_2-h_3)p_{h_1h_2h_3}f_{h_1h_3}
\\
&=& \sum_{i=1}^6 G_i.
\end{eqnarray*}
\subsubsection*{Proof that $G_3$ and $G_6$ are negligible} Using
mixing and Lemma \ref{lemaa}, we have
as $\nto$,
\begin{eqnarray*}
|G_3| &= &8 \frac{m_n^2}{n^2} \biggl|\sum_{B_3}
(n-h_1-h_2) \bigl((p_{h_1h_2}-p_0p_{h_2})+p_0
\bigl(p_{h_2}-p_0^2\bigr)+p_0^3
\bigr) f_{h_1h_2} \biggr|
\\
& \leq& cm_n^2 \sum_{h_1=r_n+1}^n
\xi_{h_1} + c \frac{m_n}{n} \sum_{h_2=r_n+1}^n
\xi_{h_2} +c \frac{(m_np_0)^3}{m_n(\sin(\lambda
_j/2))^2}=G_3'\to0.
\end{eqnarray*}
We also have
\begin{eqnarray*}
|G_6| &\le& c \frac{m_n^2}{n^2} \Biggl| \sum_{h_1=r_n+1}^{n-r_n-3}
\sum_{h_3=r_n+1}^{n-h_1-2} \Biggl( \sum
_{h_2=1}^{r_n} +\sum_{h_2=r_n+1}^{n-h_1-h_3-1}
\Biggr) (n-h_1-h_2-h_3) p_{h_1h_2h_3}
f_{h_1h_2} \Biggl|\\
&=& G_{61}+G_{62}.
\end{eqnarray*}
By \eqref{eqhelp4}, $G_{62}$ is negligible and the same arguments as
for $G_3$ show that $G_{61}\le r_n G_3'\to0$.
Thus, $G_6$ is negligible as $\nto$.
\subsubsection*{The non-negligible contributions of $G_1,G_2,G_4,G_5$.}
First, observe that
\begin{eqnarray*}
E \wh f_{nA}^2(\la_{j_1}) &= &
(m_np_0)^2 + 4m_n^2p_0
\sum_{h=1}^{r_n} \frac{n-h}{n}p_h
\cos(\lambda_{j}h)
\\
&&{}+4 \frac
{m_n^2}{n^2} \sum_{h_1=1}^{r_n}
\sum_{h_2=1}^{r_n} (n-h_1)
(n-h_2)p_{h_1}p_{h_2}f_{h_1h_2}\\
&=&
P_1+P_2+P_3 ,
\end{eqnarray*}
and we also know that \eqref{eq101} holds. Thus, \eqref{eqhe2} is
proved if we can show
that $G_1-P_1$, $G_2-P_2$ and $G_4+G_5-P_3$ are negligible.
Observe that $\cos^2\lambda=0.5(1+\cos(2\lambda))$. Then by mixing
and Lemma~\ref{lemaa},
\begin{eqnarray*}
|G_1-P_1| &= & 4 \frac{m_n^2}{n^2} \Biggl| \sum
_{h=r_n+1}^{n-1}(n-h) \bigl(\bigl(p_h-p_0^2
\bigr)+p_0^2 \bigr) 0.5 \bigl(1+ \cos (2
\lambda_{j}h)\bigr)-(m_np_0)^2 \Biggr|
\\
&\leq& c\frac{m_n}{n} m_n \sum_{h=r_n+1}^{n-1}
\xi_h +c (m_np_0)^2\Biggl |
\frac{1}{2n^2} \sum_{h=r_n+1}^{n-1}(n-h) -1 \Biggr|
+ c \frac{1}{n}(m_np_0)^2\to0.
\end{eqnarray*}

As to $G_2$, we split the index set $B_2$ into the disjoint parts for
$h_2\leq r_n$ and $h_2>r_n$.
The sum over $B_2$ restricted to $h_2>r_n$ can be shown to be bounded
by $cG_3'$.
Recall that $2f_{h_1+h_2,h_1}=\cos(\lambda_{j}h_2)+\cos(\lambda
_{j}(2h_1+h_2))$. Then
\begin{eqnarray*}
|G_2-P_2|&\leq& c G_3^{\prime}
+ \Biggl|2\frac{m_n^2}{n^2} \sum_{h_2=1}^{r_n}
\sum_{h_1=r_n+1}^{n-h_2-1}(n-h_1-h_2)
(p_{h_1h_2}+p_{h_2h_1}) \\
&&\hspace*{100pt}{}\times\bigl(\cos(\lambda_{j}h_1)+
\cos\bigl(\lambda _{j}(2h_2+h_1)\bigr)
\bigr)
\\
& &\hspace*{30pt}{} -4m_n^2p_0\sum
_{h=1}^{r_n}\frac{n-h}{n}p_h\cos(
\lambda _{j}h) \Biggr|
\\
&\leq&c G_3^{\prime}+ c\frac{r_n^2}{n}(m_np_0)^2
+ c \frac{m_nr_n}{n} m_n\sum_{h_2=r_n+1}^{n-1}
\xi_{h_2}+ c(m_np_0)^2
\frac{r_n}{n(\sin
(\lambda_{j}))^2}\to0.
\end{eqnarray*}
Here we used \eqref{eqhelp2} to rewrite $p_{h_1h_2}$ such that the mixing
condition and Lemma \ref{lemaa} can be applied.

Finally, we turn to $G_4$ and $G_5$.
By virtue of \eqref{eqhelp4} and \eqref{eqhelp5}, we can neglect
those parts of $G_4+G_5$ which contain
$(h_1,h_2,h_3)$-indices with $h_1,h_2,h_3>r_n$. Those parts of
$G_4+G_5$ for which
two indices out of $(h_1,h_2,h_3)$ exceed $r_n$ we can deal with like
$J_{23}$, and a similar argument applies
when either $h_1>r_n$ or $h_3>r_n$. Thus, we need to study those
summands in $G_4+G_5$ indexed on
$\{1\leq h_1,h_3 \leq r_n, r_n<h_2<n-h_1-h_3 \}$. We write $G_{4+5}$
for the remaining sum.
Recall that
\[
f_{h_1+h_2+h_3,h_2}+f_{h_1+h_2,h_2+h_3} = f_{h_1h_3}+\cos\bigl(\lambda
_{j}(h_1+2h_2+h_3)\bigr).
\]
Then we have
\begin{eqnarray*}
|G_{4+5}-P_3|&=& \Biggl|4 \frac{m_n^2}{n^2} \sum
_{h_1=1}^{r_n}\sum_{h_3=1}^{r_n}
\sum_{h_2=r_n+1}^{n-h_1-h_3-1}(n-h_1-h_2-h_3)
\bigl((p_{h_1h_2h_3}-p_{h_1}p_{h_3})+p_{h_1}p_{h_3}
\bigr)
\\
&&\hspace*{103pt}{}\times \bigl(2f_{h_1h_3}+2\cos\bigl(\lambda_{j_1}(h_1+2h_2+h_3)
\bigr) \bigr)
\\
&&{} -4 \frac{m_{n}^2}{n^2}\sum_{h_1=r_n+1}^{n-1}
\sum_{h_2=r_n+1}^{n-1}(n-h_1)
(n-h_2)p_{h_1}p_{h_2}f_{h_1h_2}\Biggr |
\\
&\leq& c \frac{r_n^3}{n}(m_np_0)^2+c
\frac{m_nr_n^2}{n} m_n \sum_{h_3=r_n+1}^{n-1}
\xi_{h_3}+c\frac{r_n^2}{n(\sin(\lambda
_{j_1}))^2}(m_np_0)^2
.
\end{eqnarray*}
Thus we also proved that $G_4+G_5-P_3$ is negligible.

Collecting all the arguments above, we finally proved the theorem.
\end{pf*}

\section{A discussion of related results and possible
extensions}\label{sec6}
Extremogram-type quantities for time series have been introduced by various authors. Ledford and Tawn
\cite{ledfordtawn2003} discussed $\rho_{(1,\infty)}$ as a possible
measure of
extremal dependence for univariate stationary processes
with unit Fr\'echet marginals under the regular variation condition
$P(X_0>x,X_t>x)=L_t(x)x^{-1/\eta_t}$, for slowly varying $L_t$ and
$\eta
_t\in(0,1]$.
They were particularly interested in the case of
asymptotic independence when $\rho_{(1,\infty)}(t)=0$ and
$P(X_0>x,X_t>x)/[P(X>x)]^2\to1$ as $\xto$ and also suggested diagnostic
conditions in this situation. Hill~\cite{hill2011} proposed the quantities
$\lim_{\xto} [P(X_0>x,X_t>x)/[P(X>x)]^2-1]$ as alternative measure of serial
extremal dependence in the case when the extremogram vanishes.
Fasen \textit{et al.} \cite{fasenkluppelbergschlather2009} considered
lag-dependent tail dependence coefficients under regular variation conditions
on the process $(X_t)$. These coefficients can be interpreted as
special extremograms. Hill \cite{hill2009} showed a pre-asymptotic functional
central limit theorem for the sample extremogram of univariate time
series over classes of
upper quadrants.
His mixing and domain of maximum domain of attraction are not
directly comparable with strong mixing and regular variation od
stationary sequence
s but the results are similar in spirit to Theorem 3.2 in Davis and
Mikosch \cite{davismikosch2009}, where multivariate time series can be
treated but uniform convergence  over certain classes of
sets was not considered.

Recently, various articles on the spectral analysis of
indicator functions and their covariances based on a strictly stationary
time series were written; see, for example, Dette \textit{et al.} \cite{dettehallin2011}
and the references therein, Hagemann \cite{hagemann2011}, Lee and
Subba Rao \cite{leesubba2012}.
The results are similar to those of classical time series analysis.
The aforementioned papers do not deal with the spectral analysis of
serial extremal dependence. In particular,
they do not involve sequences of indicator functions of the form
$(I_{\{a_m^{-1} X_t\in A\}})$ for sets $A$ bounded away from zero. Therefore,
these papers do not need additional conditions such as regular
variation of $(X_t)$
which are typical for extreme value theory and they do not require to
consider the
normalization $m/n$ of the periodogram but use the classical $1/n$ constants.

The present paper focuses on the basic properties of the extremal
periodogram. These properties parallel the results of classical time
series analysis, but
the proofs are different because of the triangular array nature of the
stochastic processes $(I_{\{a_m^{-1} X_t\in A\}})$. In particular, the
calculation of sufficiently high moments necessary to prove central
limit theorems
becomes rather technical. The central limit theorem for the smoothed
periodogram is still
an open question.

The (smoothed) periodogram as such contains information about the
length of
random cycles between extremal events
$\{a_m^{-1} X_t\in A\}$. But it also opens the door to the methods and
procedures of classical time series analysis, including the rich theory
related to
the integrated periodogram with applications to parameter estimation
(e.g., Whittle estimation), goodness-of-fit tests, change point analysis,
detection of long-range dependence effects and other problems.
The solution to these problems is again rather technical and
will be treated in future work.

\begin{appendix}
\section{Some trigonometric sum formulas}\label{secappendix}
Equations \eqref{eqa} and~\eqref{eqb} are given in Gradshteyn and
Ryzhik \cite{gradshteynryzhik1980}, 1.341 on page~29; \eqref{eqc}~and \eqref{eqd} are 1.352 on page~31; and \eqref{eqe1} and \eqref
{eqe2} are listed as 1.353 on page~31.
For any $\la,x$ and $n\ge1$, the following identities hold
%
\begin{eqnarray}
\sum_{k=0}^{n-1} \cos( x+k \la)&=&
\frac{\cos(
x+(n-1)\la/2)\sin(n\la/2)
}{\sin
(\lambda/2)} ,\label{eqa}
\\
\sum_{k=0}^{n-1} \sin(x+ k \la)&=&
\frac{\sin(
x+(n-1)\la/2)\sin(n\la/2)}{\sin(\la/2)} ,\label{eqb}
\\
\sum_{k=1}^{n-1}k \cos( k \la) & =&
\frac{n \sin((2n-1)\lambda
/2)}{2\sin(\lambda/2)} - \frac{1-\cos n \lambda}{4(\sin(\lambda/2))^2} ,\label{eqc}
\\
\sum_{k=1}^{n-1} k \sin( k \la) & =&
\frac{\sin(n\lambda)}{4(\sin(\lambda/2))^2} - \frac{n \cos((2n-1)\la/2)}{2\sin(\lambda/2)} ,\label{eqd}
\\
\sum_{k=1}^{n-1} p^k \sin(k
\lambda) &=& \frac{p\sin(\lambda) -
p^n \sin(n\lambda) +p^{n+1}
\sin ( (n-1)\lambda )}{1-2p\cos(\lambda) +p^2} , \label {eqe1}
\\
\sum_{k=0}^{n-1} p^k \cos(k
\lambda) &=& \frac{1-p\cos(\lambda) -
p^n \cos(n \lambda) +p^{n+1}
\cos((n-1) \lambda)}{1-2p\cos(\lambda) +p^2}. \label{eqe2}
\end{eqnarray}
Using these formulas, direct calculation yields for any frequency $\la$,\vspace*{-1pt}
%
\begin{eqnarray}
\label{eqe} &&\sum_{s=1}^{n-h} \bigl[
\cos(\lambda s)\sin\bigl(\lambda (s+h)\bigr)+\cos\bigl(\lambda(s+h)\bigr)\sin(
\lambda s) \bigr]
\nonumber
\\[-10pt]
\\[-10pt]
\nonumber
&&\quad= \sum_{s=1}^{n-h} \sin(2\lambda s+\la h)
=\frac{\sin(\la n)\sin(\la(n-h+1))}{\sin\la}-\sin(\la h).\vspace*{-1pt}
\end{eqnarray}
For distinct frequencies $\la,\w,$ we then obtain\vspace*{-1pt}
\begin{eqnarray}
\label{eqf} &&\sum_{s=1}^{n-h} \bigl[
\cos(\lambda s)\sin\bigl(\omega (s+h)\bigr)+\cos\bigl(\lambda (s+h)\bigr)\sin(
\omega s) \bigr]
\nonumber
\\[-2pt]
&&\quad= 0.5 \sum_{s=1}^{n-h} \bigl[\sin\bigl((
\la+\w)s +\w h\bigr)-\sin\bigl((\la-\w)s -\w h\bigr)\bigr]
\nonumber
\\[-2pt]
&&\qquad{}+0.5 \sum_{s=1}^{n-h} \bigl[\sin\bigl((
\la+\w)s +\la h\bigr)-\sin\bigl((\la-\w)s +\la h\bigr)\bigr]
\nonumber
\\[-2pt]
&&\quad= -\sin(\w h)
\nonumber
\\[-9pt]
\\[-9pt]
\nonumber
&&\qquad{}+0.5\frac{\sin((n-h+1)(\la+\w
)/2)}{\sin((\la+\w)/2)}\\[-2pt]
&&\qquad\phantom{+}{}\times \bigl[\sin\bigl(\w h + (n-h) (\la+\w)/2\bigr)+\sin\bigl(
\la h + (n-h) (\la+\w)/2\bigr) \bigr]
\nonumber
\\[-2pt]
&&\qquad{} -0.5 \frac{\sin((n-h+1)(\la-\w)/2)}{\sin((\la-\w)/2)}\nonumber\\[-2pt]
&&\qquad\phantom{-}{}\times \bigl[\sin \bigl(-\w h+ (n-h) (\la-\w)/2\bigr)+\sin
\bigl(\la h + (n-h) (\la-\w)/2\bigr) \bigr] ,\nonumber
\end{eqnarray}\vspace*{-23pt}
%
\begin{eqnarray}
\label{eqg} &&\sum_{s=1}^{n-h} \bigl[
\cos(\lambda s) \cos\bigl(\omega(s+h)\bigr) + \cos\bigl(\lambda(s+h)\bigr) \cos(
\omega s) \bigr]
\nonumber
\\[-2pt]
&&\quad=  0.5\sum_{s=1}^{n-h} \bigl[ \cos\bigl((
\lambda+\omega) s + \omega h\bigr) +\cos\bigl((\lambda-\omega)s - \omega h\bigr) +
\cos\bigl((\lambda+\omega)s +\lambda h\bigr)\nonumber \\[-2pt]
&&\hspace*{51pt}{}+\cos\bigl((\lambda-\omega)s +
\lambda h\bigr) \bigr]
\nonumber
\\[-2pt]
&&\quad=  -
\cos(\omega h) -\cos(\lambda h)
\nonumber
\\[-9pt]
\\[-9pt]
\nonumber
&& \qquad{}+0.5 \frac{\sin((n-h+1)(\lambda+\omega)/2)}{\sin((\lambda+
\omega)/2)} \\[-2pt]
&&\qquad\phantom{+}{}\times\bigl[ \cos\bigl(\omega h + (n-h) (\lambda+\omega)/2
\bigr) +\cos \bigl(\lambda h + (n-h) (\lambda+\omega)/2\bigr) \bigr]
\nonumber
\\[-2pt]
&&\qquad{} + 0.5 \frac{\sin((n-h+1)(\lambda-\omega)/2)}{\sin((\lambda-
\omega)/2)} \nonumber\\[-2pt]
&&\qquad\phantom{+}{}\times\bigl[ \cos\bigl(-\omega h + (n-h) (\lambda-
\omega)/2\bigr) +\cos \bigl(\lambda h + (n-h) (\lambda-\omega)/2\bigr) \bigr]
,\nonumber
\end{eqnarray}\vspace*{-23pt}
%
\begin{eqnarray}
\label{eqh} &&\sum_{s=1}^{n-h} \bigl[
\sin(\lambda s) \sin\bigl(\omega(s+h)\bigr) + \sin\bigl(\lambda(s+h)\bigr) \sin(
\omega s) \bigr]
\nonumber
\\
&&\quad =  0.5\sum_{s=1}^{n-h} \bigl[ \cos
\bigl((\lambda+\omega) s + \omega h\bigr) -\cos\bigl((\lambda-\omega)s - \omega h
\bigr) + \cos\bigl((\lambda+\omega)s +\lambda h\bigr)\nonumber \\
&&\hspace*{51pt}{}-\cos\bigl((\lambda-\omega)s
+\lambda h\bigr) \bigr]
\nonumber
\\
&&\quad=  0.5
\frac{\sin((n-h+1)(\lambda-\omega)/2)}{\sin((\lambda-
\omega)/2)} \\
&&\qquad{}\times\bigl[ \cos\bigl(-\omega h + (n-h) (\lambda-\omega)/2\bigr) +
\cos \bigl(\lambda h + (n-h) (\lambda-\omega)/2\bigr) \bigr]
\nonumber
\\
&&\qquad{} -0.5 \frac{\sin((n-h+1)(\lambda+\omega)/2)}{\sin((\lambda+
\omega)/2)}\nonumber\\
&&\qquad\phantom{-}{}\times \bigl[ \cos\bigl(\omega h + (n-h) (\lambda+\omega)/2
\bigr) +\cos \bigl(\lambda h + (n-h) (\lambda+\omega)/2\bigr) \bigr]
.\nonumber
\end{eqnarray}
Next, assume the conditions of Theorem~\ref{thmsmoth}. Then
a direct application of \eqref{eqa}--\eqref{eqd} yields for
$\lambda\in(0, \uppi)$ the following relations:
\begin{eqnarray*}
&& \sum_{s=r_n+1}^n (n-s) \sin(\lambda s+ x)
\\[-2pt]
&&\quad=  n \biggl( \frac{\sin(x+(n-1))\lambda/2 \sin(n\lambda/2)
}{\sin(\lambda/2)}\\[-2pt]
&&\hspace*{11pt}\qquad{} -\frac{\sin(x+r_n \lambda/2) \sin((r_n
+1)\lambda/2) }{\sin(\lambda/2)} \biggr)
\\[-2pt]
&&\qquad{}+ \sin(x) \biggl( \frac{n\sin((2n-1)\lambda/2)}{2\sin(\lambda
/2)} -\frac{(r_n+1)\sin((2r_n-1)\lambda/2)}{2\sin(\lambda/2)}\\[-2pt]
&&\hspace*{43pt}\qquad{} - \frac{\cos((r_n+1)\lambda)-\cos(n\lambda)}{4(\sin(\lambda
/2))^2}
\biggr)
\\[-2pt]
&&\qquad{} +\cos(x) \biggl( \frac{n\cos((2n-1)\lambda/2)}{2\sin(\lambda
/2)} -\frac{(r_n+1)\cos((2r_n-1)\lambda/2)}{2\sin(\lambda/2)}\\[-2pt]
&&\hspace*{43pt}\qquad{} - \frac{\sin(n \lambda)-\sin((r_n+1)\lambda)}{4(\sin(\lambda
/2))^2}
\biggr)
\end{eqnarray*}\vspace*{-28pt}
\begin{eqnarray*}
&& \sum_{s=r_n+1}^n (n-s) \cos(\lambda s+ x)
\\
&&\quad=  n \biggl( \frac{\cos(x+(n-1)\lambda/2) \sin(n\lambda/2)
}{\sin(\lambda/2)} -\frac{\cos(x+r_n \lambda/2) \sin((r_n
+1)\lambda/2) }{\sin(\lambda/2)} \biggr)
\\
&&\qquad{}- \cos(x) \biggl( \frac{n\sin((2n-1)\lambda/2)}{2\sin(\lambda
/2)} -\frac{(r_n+1)\sin((2r_n-1)\lambda/2)}{2\sin(\lambda/2)} \\
&&\hspace*{43pt}\qquad{}- \frac{\cos((r_n+1)\lambda)-\cos(n\lambda)}{4(\sin(\lambda
/2))^2}
\biggr)
\\
&&\qquad{} +\sin(x) \biggl( \frac{n\cos((2n-1)\lambda/2)}{2\sin(\lambda
/2)} -\frac{(r_n+1)\cos((2r_n-1)\lambda/2)}{2\sin(\lambda/2)}\\
&&\hspace*{43pt}\qquad{} - \frac{\sin(n \lambda)-\sin((r_n+1)\lambda)}{4(\sin(\lambda
/2))^2}
\biggr).
\end{eqnarray*}

\begin{lemma}\label{lemaa}
Under the assumptions of Theorem~\ref{thmsmoth}
the following relations hold uniformly for
$\lambda\in(0,2\uppi)$, as $\nto$,
\begin{eqnarray*}
&&\sum_{h=r_n+1}^{n-1} (n-h) \cos(
\lambda h+ x)
\\
&&\quad= \frac{n\cos(x+(n-1)\lambda/2)\sin(n\lambda/2)}{\sin(\lambda
/2)}-n-\frac{n\sin((2n-1)\lambda/2)}{2\sin(\lambda/2)}+\frac
{1-\cos(n\lambda)}{4(\sin(\lambda/2))^2}
\\
& &\qquad{}- \frac{n\cos(x+r_n\lambda/2)\sin((r_n+1)\lambda/2)}{\sin
(\lambda/2)}-n+\frac{(r_n+1)\sin((2r_n+1)\lambda/2)}{2\sin(\lambda
/2)}\\
&&\qquad{}-\frac{1-\cos(r_n+1)\lambda}{4(\sin(\lambda/2))^2}
\\
&&\quad=  \mathrm{O}\bigl(n/\bigl(\sin(\la/2)\bigr)^2\bigr) ,
\\
&& \sum_{h=r_n+1}^{n-1} (n-h) \sin(
\lambda h+ x)
\\
&&\quad= \frac{n\sin(x+(n-1)\lambda/2)\sin(n\lambda/2)}{\sin(\lambda
/2)}-\frac{\sin(n\lambda)}{4(\sin(\lambda/2))^2}+\frac{n\cos
((2n-1)\lambda/2)}{2\sin(\lambda/2)}
\\
& &\qquad{}-\frac{n\sin(x+r_n\lambda/2)\sin((r_n+1)\lambda/2)}{\sin
(\lambda/2)}+\frac{\sin(r_n\lambda)}{4(\sin(\lambda/2))^2}-\frac
{n\cos((2r_n+1)\lambda/2)}{2\sin(\lambda/2)}
\\
&&\quad=  \mathrm{O}\bigl(n/\bigl(\sin(\la/2)\bigr)^2\bigr).
\end{eqnarray*}
\end{lemma}
%
\section{The spectral density $f_A$ of an $\operatorname{ARMA}(1,1)$ process}\label{exam1}

In this section, we calculate the spectral density
$f_A$ for an $\operatorname{ARMA}(1,1)$ process and the set $A=(1,\infty)$. The process
$(X_t)$ is given as
the stationary causal solution to the difference equation
\[
X_t = \phi X_{t-1} +Z_t + \theta
Z_{t-1} , \qquad t\in\bbz,
\]
where $0<|\phi|<1$ and $\theta\in\bbr$. From Brockwell and Davis
\cite{brockwelldavis2002}, (2.3.3), we obtain the coefficients
$(\psi_j)$ of the linear process representation of $(X_t)$ (cf.
\eqref{eq8}):
\[
\psi_0 =1 , \qquad\psi_j = \phi^{j-1} (\phi+
\theta) ,\qquad j\geq1.
\]
We assume that $(Z_t)$ is an i.i.d. regularly varying sequence with index $\alpha>0$.\\[1mm]
\emph{The case $\phi\in(0,1)$, $\theta+\phi>0$, $p>0$.} A direct
application
of \eqref{eq10} yields that
\begin{eqnarray*}
\rho_A(h)&=& \frac{\min(1,
\psi_{h}^\alpha) +\sum_{i=h+1}^\infty\psi_i^\alpha} {
\sum_{i=0}^ \infty\psi_i^\alpha}
\\
&=&\frac{\min(1,\phi^{\alpha(h-1)}(\theta+\phi)^\alpha)
+\phi^{\alpha h} (\theta+\phi)^\alpha(1-\phi^\alpha)^{-1}} {
1+ (\theta+\phi)^\alpha(1-\phi^\alpha)^{-1}} ,\qquad h\ge1.
\end{eqnarray*}
Define
$h_0=\min\{h\ge0\dvt \phi^{\alpha h}(\theta+\phi)^\alpha< 1\}$
and write (see Appendix~\ref{secappendix})
\begin{eqnarray*}
L^{(1)} (n,x,\lambda)&=& \sum_{h=1}^{n}
\cos(x+h\lambda)\\
&=& \cases{ \displaystyle\frac{\cos(x+n \lambda)\sin((n +1)\lambda/2)}{\sin(\lambda/2)} -1 ,&\quad $n\ge1,$ \vspace*{2pt}
\cr
0 ,
& \quad $n=0$;}
\end{eqnarray*}
\begin{eqnarray*}
&&L^{(2)} (n,x,\alpha,\lambda) \\
&&\quad= \sum_{h=1}^{n}
|\phi|^{\alpha h} \cos(x+h\lambda)
\\
&&\quad = \cases{\displaystyle\frac{|\phi|^\alpha\cos(x+\lambda) - |\phi|^{2\alpha} \cos
(x)-|\phi|^{\alpha(n+1)} \cos(x+(n+1) \lambda) +
|\phi|^{\alpha(n+2)}\cos(x+n \lambda)}{|1-|\phi|^\alpha\ex
^{-\mathrm{i}\lambda}|^2} ,\vspace*{2pt}\cr
\qquad \hspace*{10.5pt}n\ge1 , \vspace*{2pt}
\cr
0 , \qquad  n=0,
\vspace*{2pt}
\cr
\displaystyle \frac{|\phi|^{\alpha} \cos(x+\lambda)-|\phi|^{2\alpha}
\cos(x)}{|1-|\phi|^\alpha\ex^{-\mathrm{i}\lambda}|^2} ,\vspace*{2pt}\cr
\qquad\hspace*{10.5pt}  n=\infty. }
\end{eqnarray*}
Then
\[
\rho_A (h) = \cases{ %
c^{(1)}_{\alpha} (\phi, \theta) +\phi^{\alpha
h}
c_\alpha^{(2)}(\phi,\theta) , &\quad  $h\le h_0 ,$
\vspace*{2pt}\cr
\phi^{\alpha(h-1)}c_\alpha^{(2)}(\phi,\theta)
,&\quad $h>h_0 ,$ }
\]
where
\begin{eqnarray*}
c^{(1)}_{\alpha} (\phi, \theta) = \frac{1-\phi^\alpha}{1-\phi
^\alpha+(\phi+\theta)^\alpha} \quad\mbox{and}\quad
c^{(2)}_{\alpha} (\phi, \theta) = \frac{(\phi+\theta)^\alpha
}{1-\phi^\alpha+
(\phi+\theta)^\alpha}.
\end{eqnarray*}
The corresponding spectral density is given by
\begin{eqnarray*}
f_A(\lambda) &=& 1+ 2c^{(1)}_\alpha(\phi,
\theta) \sum_{h=1}^{h_0} \cos(h\lambda) + 2
\bigl(1-\phi^{-\alpha
}\bigr)c^{(2)}_\alpha(\phi, \theta)
\sum_{h=1}^{h_0} \phi^{\alpha h} \cos(h
\lambda)
\\
&&{} + 2\phi^{-\alpha} c^{(2)}_\alpha(\phi, \theta) \sum
_{h=1}^{\infty} \phi^{\alpha h} \cos(h
\lambda)
\\
&=& 1+2 c^{(1)}_\alpha(\phi, \theta) L^{(1)}
(h_0,0, \lambda) + 2 \bigl(1-\phi^ {-\alpha}\bigr)
c^{(2)}_\alpha(\phi, \theta) L^{(2)}
(h_0,0,\alpha, \lambda)
\\
&&{}+ 2\phi^{-\alpha}c^{(2)}_\alpha(\phi,
\theta)L^{(2)} (\infty,0,\alpha,\lambda).
\end{eqnarray*}
\emph{The case $\phi\in(0,1)$, $\theta+\phi<0$, $q>0$.}
In view of \eqref{eq10}, we have
\begin{eqnarray*}
\rho_A (h) &=& \frac{q\sum_{i=0}^{\infty} \phi^{\alpha h +\alpha
i} |\phi+\theta|^\alpha}{p+ q\sum_{i=0}^{\infty} \phi^{\alpha i}
|\phi+\theta|^\alpha} = \frac{q \phi^{\alpha h} |\phi+\theta
|^\alpha}{p(1-\phi^\alpha) + q|\phi+\theta|^\alpha} = \phi
^{\alpha h} c^{(3)}_\alpha(\phi, \theta),\qquad  h\ge1 ,
\\
f_A(\lambda) &=& 1+ 2 c^{(3)}_\alpha(\phi,
\theta) L^{(2)} (\infty , 0,\alpha,\lambda).
\end{eqnarray*}
\emph{The case $\phi\in(-1,0)$, $\theta+\phi>0 $, $p>0$.}
If $h=2k+1$ for integer $k\ge0$ the summand $p (\min(\psi_{i}^+ ,
\psi_{i+h}^+))^\alpha+ q (\min(\psi_{i}^- ,
\psi_{i+h}^-))^\alpha$ in \eqref{eq10} vanishes for $i\ge1$. Thus,
\[
\rho_{A}(h) = \frac{p\min(1,
|\psi_{h}|^\alpha)} {
p + \sum_{i=1}^{\infty}  [ p |\psi_{2i-1}|^\alpha+ q |\psi
_{2i}|^\alpha ] }.
\]
For $h=2k>0$,
\[
\rho_{A} (h) = \frac{\sum_{i=1}^{\infty}  [ p |\psi
_{2i+h-1}|^\alpha+ q |\psi_{2i+h}|^\alpha ] }{ p + \sum_{i=1}^{\infty}  [ p |\psi_{2i-1}|^\alpha+ q |\psi_{2i}|^\alpha
] }.
\]
Define $
k_1=\min\{k\ge0\dvt |\phi|^{2 k}(\theta+\phi) <1\}$.
Then,
\[
\rho_A (h)= \cases{ %
c^{(4)}_{\alpha} (\phi, \theta) , & \quad $h=2k-1 , 1\le k\le
k_1 ,$
\vspace*{2pt}\cr
\phi^{\alpha(h-1)} c_\alpha^{(5)}(\phi,\theta) ,&\quad $h=2k-1 ,
k>k_1 ,$
\vspace*{2pt}\cr
\phi^{\alpha h} c_\alpha^{(6)}(\phi,\theta) ,&\quad $h=2k , k
\ge1 ,$}
\]
where
\begin{eqnarray*}
c^{(4)}_{\alpha} &=&\frac{p(1-|\phi|^{2\alpha})}{p(1-|\phi
|^{2\alpha} +(\phi+\theta)^{\alpha}) +q|\phi|^{\alpha} (\phi
+\theta)^{\alpha}} ,
\\
 c^{(5)}_{\alpha}&=&\frac{p(\phi+\theta)^{\alpha} (1-|\phi
|^{2\alpha})}{p(1-|\phi|^{2\alpha} +(\phi+\theta)^{\alpha})
+q|\phi|^{\alpha} (\phi+\theta)^{\alpha}} ,
\\
 c^{(6)}_{\alpha} &=& \frac{p(\phi+\theta)^{\alpha} +q|\phi
|^\alpha(\phi+\theta)^{\alpha}}{p(1-|\phi|^{2\alpha} +(\phi
+\theta)^{\alpha}) +q|\phi|^{\alpha} (\phi+\theta)^{\alpha}}.
\end{eqnarray*}
The corresponding spectral density is
\begin{eqnarray*}
f_{A}(\lambda) &=& 1+2c^{(4)}_{\alpha} (\phi,
\theta) \sum_{k=1}^{k_1} \cos\bigl((2k-1)
\lambda\bigr) + 2|\phi|^{-2\alpha
}c^{(5)}_{\alpha} (\phi,
\theta) \sum_{k=k_1 +1}^{\infty} |\phi
|^{\alpha(2k)} \cos\bigl((2k-1)\lambda\bigr)
\\
&&{} + 2c^{(6)}_{\alpha} (\phi, \theta) \sum
_{k=1}^{\infty} |\phi|^{2k\alpha} \cos(2k\lambda)
\\
&=& 1+ 2c^{(4)}_{\alpha} (\phi, \theta) L^{(1)}
(k_1,-\lambda,2 \lambda) \\
&&{}+2|\phi|^{-2\alpha}c^{(5)}_{\alpha}
\bigl[ L^{(2)} (\infty , -\lambda, 2\alpha, 2\lambda) -L^{(2)}
(k_1 , -\lambda, 2\alpha, 2\lambda) \bigr]
\\
&& {}+ 2c^{(6)}_{\alpha} (\phi, \theta) L^{(2)} (\infty,
0 , \alpha, 2\lambda).
\end{eqnarray*}
\emph{The case $\phi\in(-1,0)$, $\theta+\phi<0 $, $p>0$.} If
$h=2k+1$ for integer $k\ge0$
the summand $p (\min(\psi_{i}^+ ,
\psi_{i+h}^+))^\alpha+ q (\min(\psi_{i}^- ,
\psi_{i+h}^-))^\alpha$ in \eqref{eq10} vanishes for $i\ge0$. Thus,
\[
\rho_{A} (h) =0.
\]
For $h=2k >0$,
\[
\rho_{A}(h) = \frac{p\min(1, |\psi_{h}|^\alpha)+\sum_{i=0}^{\infty}  [ p|\psi_{2i+h+2}|^\alpha+q|\psi
_{2i+h+1}|^\alpha ] }{\sum_{i=0}^{\infty}  [ p|\psi
_{2i}|^\alpha+q|\psi_{2i+1}|^\alpha ]}.
\]
Define $
k_2=\min\{k\ge0\dvt |\phi|^{2 k+1}| \theta+\phi| <1\}$. Then
\begin{eqnarray*}
\rho_{A}(2k)=\cases{ c_\alpha^{(7)} +|
\phi|^{2\alpha k} c_\alpha^{(8)} , & \quad $k\le k_2$
,\vspace*{2pt}
\cr
|\phi|^{2\alpha k} c_\alpha^{(9)} , &\quad
$k> k_2$ ,}
\end{eqnarray*}
where
\begin{eqnarray*}
 c_\alpha^{(7)} &=& \frac{p(1-|\phi|^{2\alpha})}{p(1-|\phi
|^{2\alpha} ) +p|\phi|^\alpha|\phi+\theta|^\alpha+q |\phi+\theta
|^\alpha} ,
\\
 c_\alpha^{(8)} &=& \frac{p|\phi|^\alpha|\phi+\theta|^\alpha+q
|\phi+\theta|^\alpha}{p(1-|\phi|^{2\alpha} ) +p|\phi|^\alpha
|\phi+\theta|^\alpha+q |\phi+\theta|^\alpha} ,
\\
 c_\alpha^{(9)} &=& \frac{p|\phi|^{-\alpha} |\phi+\theta|^\alpha
+q |\phi+\theta|^\alpha}{p(1-|\phi|^{2\alpha} ) +p|\phi|^\alpha
|\phi+\theta|^\alpha+q |\phi+\theta|^\alpha}.
\end{eqnarray*}
The corresponding spectral density is
\begin{eqnarray*}
f_{A}(\lambda) &=& 1+2 c_\alpha^{(7)} \sum
_{k=1}^{k_2} \cos (2k\lambda) +2
\bigl(c_\alpha^{(8)} -c_\alpha^{(9)} \bigr)
\sum_{k=1}^{k_2} |\phi|^{2k\alpha}
\cos(2k\lambda) +2c_\alpha^{(9)} \sum
_{k=1}^{\infty} |\phi|^{2k\alpha} \cos(2k\lambda)
\\
&=& 1+2c_\alpha^{(7)} L^{(1)} (k_2,0,2
\lambda) +2\bigl(c_\alpha^{(8)} -c_\alpha^{(9)}
\bigr)L^{(2)} (k_2,0,2\alpha, 2\lambda)+ 2 c_\alpha
^{(9)} L^{(2)} (\infty, 0,2\alpha, 2\lambda).
\end{eqnarray*}%
\end{appendix}%
%
%
\section*{Acknowledgments}
We would like to thank the reviewers of our paper
for careful reading and comments, in particular for pointing out
several useful references. Thomas Mikosch's research is partly supported
by the Danish Research Council (FNU) Grants
09-072331 ``Point process modelling and statistical inference''
and 10-084172 ``Heavy tail phenomena: Modeling and estimation''.
The research of Yuwei Zhao is supported by the Danish Research Council
Grant 10-084172.
%

\printhistory
\end{document}